\documentclass{amsart}

\usepackage{amsmath}
\usepackage{amssymb}
\usepackage{amsthm}
\usepackage{mathrsfs}
\usepackage{pdfpages}
\usepackage{dsfont}
\usepackage{lineno}
\usepackage{enumerate}
\usepackage{graphicx}
\usepackage{epstopdf}
\usepackage{color}
\usepackage{xspace}
\usepackage{subfig}
\usepackage{pdfsync}

\DeclareMathOperator*{\argmin}{argmin}

\DeclareGraphicsExtensions{.pdf}
\graphicspath{{./PDF/}}

\newcommand{\bq}{\begin{equation}}
\newcommand{\eq}{\end{equation}}
\newcommand{\R}{\mathbb{R}}

\newcommand{\abs}[1]{\left\vert#1\right\vert}

\newcommand{\bO}{\mathcal{O}}
\newcommand{\bF}{\mathcal{F}}

\newcommand{\Dt}{\mathcal{D}}

\newcommand{\QW}{quadratic Wasserstein metric\xspace}
\newcommand{\MA}{Monge-Amp\`ere\xspace}

\newcommand{\diag}{\text{diag}}
\newcommand{\vp}{v^\perp}
\newcommand{\M}{\mathcal{M}}
\newcommand*\Laplace{\mathop{}\!\mathbin\bigtriangleup}

\newtheorem{theorem}{Theorem}
\theoremstyle{lemma}

\newtheorem{definition}{Definition}
\newtheorem{proposition}[theorem]{Proposition}
\theoremstyle{remark}

\makeindex

\begin{document}
\title{Seismic Imaging and Optimal Transport}

\author{Bj\"orn Engquist}
\address{Department of Mathematics and ICES, The University of Texas at Austin, 1 University Station C1200, Austin, TX 78712 USA}
\email{engquist@math.utexas.edu}

\author{Yunan Yang}
\address{Department of Mathematics, The University of Texas at Austin, 1 University Station C1200, Austin, TX 78712 USA}
\email{yunanyang@math.utexas.edu}

\begin{abstract}
Seismology has been an active science for a long time. It changed character about 50 years ago when the earth’s vibrations could be measured on the surface more accurately and more frequently in space and time. The full wave field could be determined, and partial differential equations (PDE) started to be used in the inverse process of finding properties of the interior of the earth. We will briefly review earlier techniques but mainly focus on Full Waveform Inversion (FWI) for the acoustic formulation. FWI is a PDE constrained optimization in which the variable velocity in a forward wave equation is adjusted such that the solution matches measured data on the surface. The minimization of the mismatch is usually coupled with the adjoint state method, which also includes the solution to an adjoint wave equation. The least-squares ($L^2$) norm is the conventional objective function measuring the difference between simulated and measured data, but it often results in the minimization trapped in local minima. One way to mitigate this is by selecting another misfit function with better convexity properties. Here we propose using the quadratic Wasserstein metric ($W_2$) as a new misfit function in FWI. The optimal map defining $W_2$ can be computed by solving a Monge-Amp\`ere equation. Theorems pointing to the advantages of using optimal transport over $L^2$ norm will be discussed, and a number of large-scale computational examples will be presented. 
\end{abstract}

\date{\today}
\maketitle

\noindent \textbf{Keywords.} Seismic Imaging, Full-waveform Inversion, Optimal Transport, \MA equation

\noindent \textbf{Math Subject Classification.} 65K10, 65Z05, 86A15,  86A22 

\tableofcontents

\section{Introduction}
\begin{figure}[tbp]
    \subfloat[]{\includegraphics[height=0.4\textwidth]{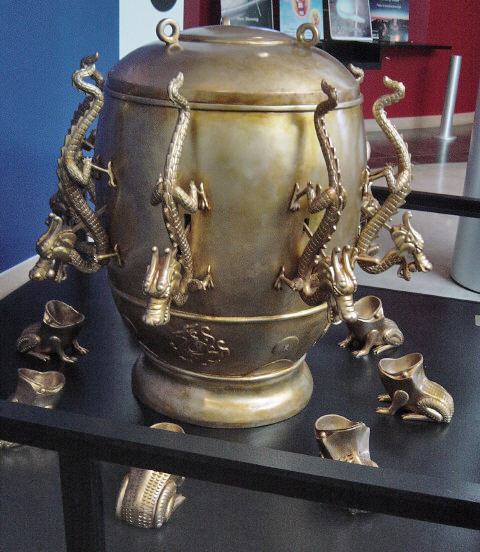}
\label{fig:scope} }
\hspace{.2cm}
    \subfloat[]{\includegraphics[height=0.4\textwidth]{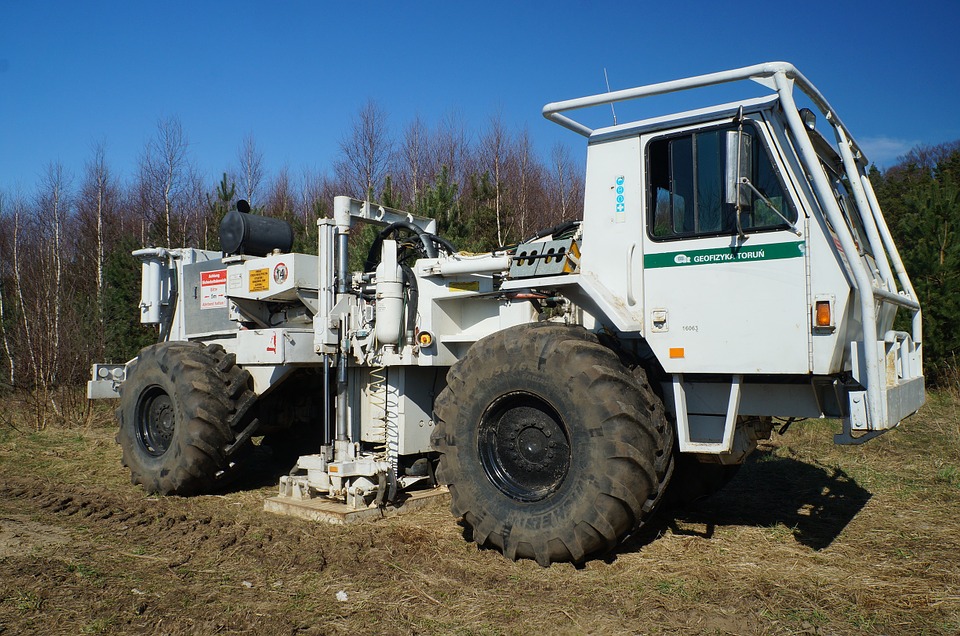}\label{fig:survey}}
\caption{(a)~The first Seismoscope designed in AD 132  and (b)~Modern seismic vibrator used in seismic survey}
\end{figure}

Earth Science is an early scientific subject. The efforts started as early as AD 132 in China when Heng Zhang invented the first seismoscope in the world (Figure~\ref{fig:scope}). The goal was to record that an earthquake had happened and to try to determine the direction of the earthquake. Substantial progress in seismology had to wait until about 150 years ago when seismological instruments started to record travel time.

With increasing sophistication in devices measuring the vibrations of seismic waves and in the availability of high-performance computing increasingly advances mathematical techniques could be used to explore the interior of the earth. The development started with calculations by hand based on geometrical optics and travel time measurement. It continued with a variety of wave equations when the equipment allowed for measuring wave fields and modern computers became available. As we will see below a wide range of mathematical tools are used today in seismic imaging, including partial differential equation (PDE) constrained optimization, advanced signal processing, optimal transport and the \MA equation.

Since 19th-century modern seismographs were developed to record seismic signals, which are vibrations in the earth. In 1798 Henry Cavendish measured the density of the earth with less than 1\% error compared with the number we can measure nowadays. Nearly one hundred years later, German physicist Emil Wiechert first discovered that the earth has a layered structure and the theory was further completed as the crust-mantle-core three-layer model in 1914 by one of his student Beno Gutenberg. In the meantime, people studied the waves including body waves and surface waves to better understand the earthquake. P-waves and S-waves were first clearly identified for their separate arrivals by English geologist Richard Dixon Oldham in 1897. The Murchison earthquake in 1929 inspired the Danish female seismologist and geophysicist Inge Lehmann to study the unexpected P-waves recorded by the seismographs. Later on, she proposed that the core of the earth has two parts: the solid inner core of iron and a liquid outer core of nickel-iron alloy, which was soon acknowledged by peer geophysicists worldwide. 




We will see that measuring travel time plays a vital role in the development of modern techniques for the inverse problem of finding geophysical properties from measurements of seismic waves on the surface. The methods are often related to travel time tomography. They are quite robust and cost-efficient for achieving low-resolution information of the subsurface velocities. The forward problem is based on ray theory or geometric optics~\cite{Bijwaard1998,Uhlmann2001}. 

\begin{figure}[tbp]
    \subfloat[]{\includegraphics[height=0.45\textwidth]{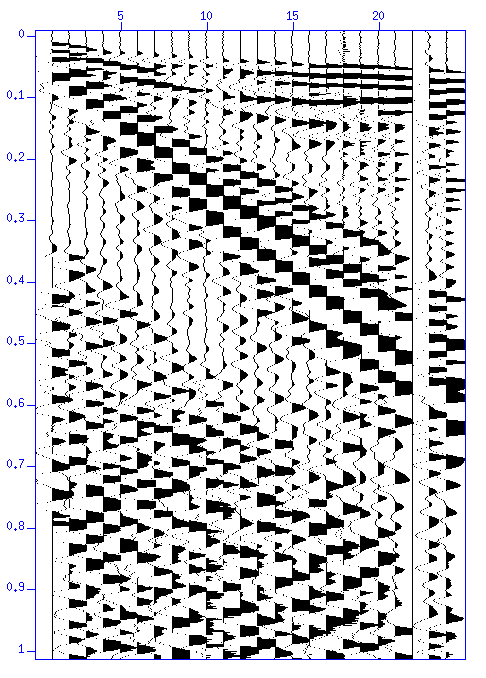}
\label{fig:data} }
\hspace{.2cm}
    \subfloat[]{\includegraphics[height=0.45\textwidth]{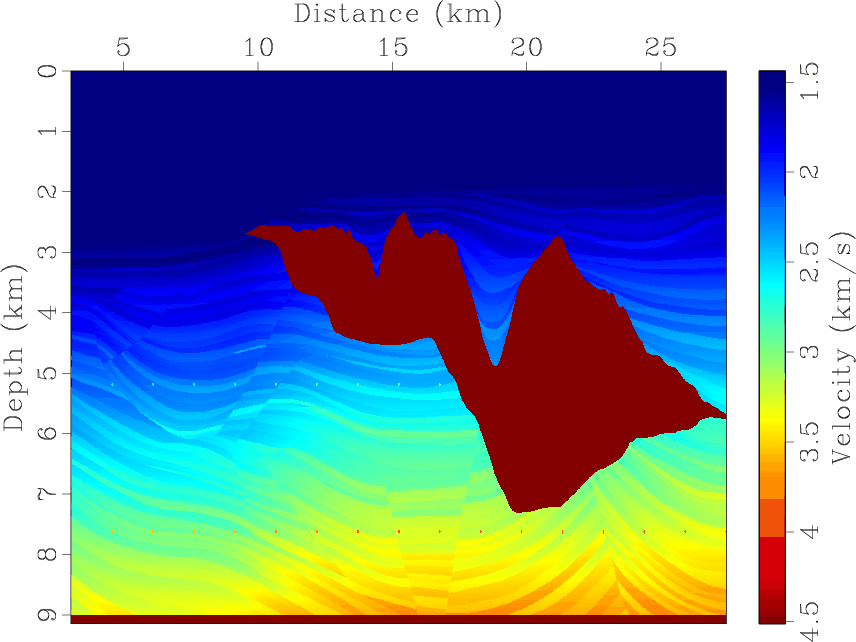}\label{fig:vel}}
\caption{(a)~An example of the seismic data measured from the receivers  and (b)~Goal of inversion: geophysical properties as in the Sigsbee velocity model~\cite{Bashkardin}}
\end{figure}

The development of man-made seismic sources and advanced recording devices (Figure~\ref{fig:survey}) facilitate the research on the entire wavefields in time and space (Figure~\ref{fig:data}) rather than merely travel time. This setup results in a more controlled setting and large amounts of data, which is needed for an accurate inverse process of estimating geophysical properties, for example, Figure~\ref{fig:vel}. The forward modeling is a wave equation with many man-made sources and many receivers.
The wave equation can vary from pure acoustic waves to anisotropic viscoelasticity. Even if there are various techniques in computational exploration seismology, there are two processes that currently stand out: reverse time migration (RTM)~\cite{Baysal1983, Yoon2004} and full waveform inversion (FWI)~\cite{tarantola1987inverse,tarantola1982generalized}.

Migration techniques can be applied in both the time domain and the frequency domain following the early breakthroughs by Claerbout on imaging conditions~\cite{Claerbout1971, Claerbout}. In RTM the computed forward wavefield starting from the source is correlated in time with the computed backward wavefield which is modeled with the measured data as the source term in the adjoint wave equation. The goal is to determine details of the reflecting surfaces as, for example, faults and sedimentary layers based on the measured data and a rough estimate of the geophysical properties. The least-squares reverse time migration (LSRTM)~\cite{Dai2012} is a new migration technique designed to improve the image quality generated by RTM. Reflectivity is regarded as a small perturbation in velocity, and the quantity is recovered through a linear inverse problem.

FWI is a high-resolution seismic imaging technique which recently gets great attention from both academia and industry~\cite{Virieux2017}. The goal of FWI is to find both the small-scale and large-scale components which describe the geophysical properties using the entire content of seismic traces. A trace is the time history measured at a receiver.
In this paper, we will consider the inverse problem of finding the wave velocity of an acoustic wave equation in the interior of a domain from knowing the Cauchy boundary data together with natural boundary conditions~\cite{Clayton1977}, which is implemented by minimizing the difference between computed and measured data on the boundary. It is thus a PDE-constrained optimization. 

There are various kinds of numerical techniques that are used in seismic inversion, but FWI is increasing in popularity even if it is still facing three main computational challenges. First, the physics of seismic waves are complex, and we need more accurate forward modeling in inversion going from pure acoustic waves to anisotropic viscoelasticity~\cite{virieux2009overview}.  
Second, even as PDE-constrained optimization, the problem is highly non-convex. FWI requires more efficient and robust optimization methods to tackle the intrinsic nonlinearity. Third, the least-squares norm, classically used in FWI, suffers from local minima trapping, the so-called cycle skipping issues, and sensitivity to noise~\cite{Seismology2011}. We will see that optimal transport based Wasserstein metric is capable of dealing with the last two limitations by including both amplitudes mismatches and travel time differences~\cite{EFWass, engquist2016optimal}.


We will introduce the mathematical formulation of these techniques in the following sections. The emphasis will be on FWI, but we will also summarize the state of the art of other standard imaging steps. Finally, we will relate FWI to RTM and LSRTM. These approaches all involve the interaction of the forward and the time-reversed wavefields, which is well known as the ``imaging condition'' in geophysics.

\section{Seismic Imaging}
Seismic data contains interpretable information about subsurface properties. Imaging predicts the spatial locations as well as specifies parameter values describing the earth properties that are useful in seismology. It is particularly important for exploration seismology which mainly focuses on prospecting for energy sources, such as oil, gas, coal. Seismic attributes contain both travel time records and waveform information to create an image of the subsurface to enable geological interpretation, and to obtain an estimate of the distribution of material properties in the underground. Usually, the problem is formulated as an inverse problem incorporating both physics and mathematics. Seismic inversion and migration are terms often used in this setting. 
\subsection{Seismic data}
There are two types of seismic signals. Natural earthquakes propagate with substantial ultra-low frequency wave energy and penetrate deeply through the whole earth. Recorded by seismometers, the natural seismic waves are used to study earth structures.
The other type of data is generated by man-made ``earthquakes'' to obtain an image of the sedimentary basins in the interior of the earth close to the surface. A wavefield has to be produced using suitable sources at appropriate locations, measured by receivers at other locations after getting reflected back from within the earth, and stored using recorders. 

In this paper, we mainly discuss the second type of seismic events.
The raw seismic data is not ideal to interpret and to create an accurate image of the subsurface. Recorded artifacts are related to the surface upon which the survey was performed, the instruments of receiving and recording and the noise generated by the procedure.
We must remove or at least minimize these artifacts. Seismic data processing aims to eliminate or reduce these effects and to leave only the influences due to the structure of geology for interpretation. Typical data processing steps include but are not limited to deconvolution, demultiple, deghosting, frequency filtering, normal moveout (NMO) correction, dip moveout (DMO) correction, common midpoint (CMP) stack, vertical seismic profiling (VSP), etc~\cite{RoyChowdhury2011, Yilmaz2001}.

In the recent two decades, the availability of the increased computer power makes it possible to process each trace of the recorded common source gathers separately, aiming for a better image. We will discuss several primary imaging methods such as traveltime tomography, seismic migration, least squares migration and full waveform inversion (FWI). 

\subsection{Traveltime tomography}
Most discoveries related to the structure of the earth were based on the assumption that seismic waves can be represented by rays, which is closely associated with geometric optics~\cite{Rawlinsona, Rawlinson, Zelt1999}. 
The primary advantages are its applicability to complex, isotropic and anisotropic, laterally varying layered media and its numerical efficiency in such computations. A critical observation is the travel time information of seismic arrivals. We can understand many arrival time observations with ray theory~\cite{cerveny1977ray}, which describes how short-wavelength seismic energy propagates. 

As a background illustration, we will derive the ray tracing expressions in a 1D setting where the velocity only varies vertically~\cite{shearer2009introduction}. Ray tracing in general 3D structure is more complicated but follows similar principles. Considering a laterally homogeneous earth model where velocity $v$ only depends on depth, the ray parameter which is also called the horizontal slowness $p$, can be expressed in the following equation by the Snell's law:
\bq
p = s(z) \sin(\theta) = \frac{dT}{dX},
\eq
where $s(z)$ ($= \frac{1}{v(z)}$) is the slowness, $\theta$ is the incidence angle, $T$ is the travel time, $X$ is the horizontal range. At the turning point depth $z_p$, $ p = s(z_p) $, a constant for a given ray. The vertical slowness $\eta = \sqrt{s^2 - p^2}$.

When the velocity is a continuous function of depth, the surface to surface travel time  $T(p)$ and the distance traveled $X(p)$ have the following expressions:
\bq\label{eq:travelT}
T(p) = 2\int_0^{z_p} \frac{s^2(z)}{\sqrt{s^2(z) - p^2}} dz = 2\int_0^{z_p} \frac{s^2(z)}{\eta} dz,
\eq
and
\bq\label{eq:travelX}
X(p) = 2p\int_0^{z_p} \frac{dz}{\sqrt{s^2(z) - p^2}} = 2p\int_0^{z_p} \frac{dz}{\eta}.
\eq

The expressions above are the forward problem in traveltime tomography. The seismologists are interested in inverting model parameter $s(z)$ from observed traveltime $T$ and traveled distance $X$. 
Using integral transform pair, we can obtain
\bq~\label{eq:travelV}
z(s) =  - \frac{1}{\pi} \int_{s_0}^{s} \frac{X(p)}{\sqrt{p^2 - s^2(z) }} d(p) = \frac{1}{\pi} \int_0^{X(s)} \cosh^{-1}(p/s) dX,
\eq
which gives us the 1D velocity model.

Equation~\eqref{eq:travelV} is one example of the 1D velocity inversion problem at a given depth. There are limitations about traveltime tomography in general. First, the first arrivals are inherently nonunique. Second, the lateral velocity variations are not considered in this setting. If we divide the earth model into blocks, the 3D velocity inversion techniques can resolve some of the lateral velocity perturbations by using the travel time in each block. The problem can be formulated into a least-squares ($L^2$) inversion by minimizing the travel time residual between the predicted time and the observed time: $ ||t_\text{obs} - t_\text{pred}||_2^2$~\cite{shearer2009introduction,Zelt2011}.

One limitation of ray theory is that it is applicable only to smooth media with smooth interfaces, in which the characteristic dimensions of inhomogeneities are considerably larger than the dominant wavelength of the considered waves. The ray method can yield distorted results and will fail at caustics or in general at so-called singular regions~\cite{Cerveny2011}. Moreover, much more information is available from the observed seismograms than travel times. To some extent, travel time tomography can be seen as phase-based inversion, and next, we will introduce waveform-based methods where the wave equation plays a significant role.

\subsection{Reverse Time Migration}
To overcome the difficulties of ray theory and further improve image resolutions,  reverse time migration (RTM),  least-squares reverse time migration (LSRTM) and full-waveform inversion (FWI) replace the semi-analytical solutions to the wave equation by fully numerical solutions including the full wavefield. Without loss of generality, we will explain all the methods in a simple acoustic setting:
\begin{equation}\label{eq:FWD}
     \left\{
     \begin{array}{rl}
     & m(\mathbf{x})\frac{\partial^2 u(\mathbf{x},t)}{\partial t^2}- \Laplace u(\mathbf{x},t) = s(\mathbf{x},t)\\
    & u(\mathbf{x}, 0 ) = 0                \\
    & \frac{\partial u}{\partial t}(\mathbf{x}, 0 ) = 0    \\
     \end{array} \right.
\end{equation}
We assume the model $m(\mathbf{x}) = \frac{1}{c(\mathbf{x})^2}$ where $c(\mathbf{x})$ is the velocity, $u(\mathbf{x},t)$ is the wavefield, $s(\mathbf{x},t)$ is the source. It is a linear PDE but a nonlinear operator from model domain $m(\mathbf{x})$ to data domain $u(\mathbf{x},t)$. 

Despite the fact that migration can be used to update velocity model~\cite{liu1995migration, sava2004wave, symes2008migration}, its chief purpose is to transform measured reflection data into an image of reflecting interfaces in the subsurface. There are two principal varieties of migration techniques: reverse time migration (RTM) which gives a modest resolution of the reflectivity~\cite{Baysal1983, yoon2004challenges} and least-squares reverse-time migration (LSRTM) which typically yields a higher resolution of the reflectivity~\cite{Dai2012,Dai2013}. 

Reverse-time migration is a prestack two-way wave-equation migration to illustrate complex structure, especially strong contrast geological interfaces such as environments involving salts. Conventional RTM uses an imaging condition which is the zero time-lag cross-correlation between the source and the receiver wavefields~\cite{Claerbout1971}:
\bq\label{eq:IC}
R(\mathbf{x}) = \sum_\text{shots}\int_0^T u(\mathbf{x},t)\cdot v(\mathbf{x},t) dt,
\eq
where $u$ is the source wavefield in~\eqref{eq:FWD} and $v$ is the receiver wavefield which is the solution to the adjoint equation~\eqref{eq:rtm_adj}:
\begin{equation} \label{eq:rtm_adj}
     \left\{
     \begin{array}{rl}
     & m(\mathbf{x})\frac{\partial^2 v(\mathbf{x},t)}{\partial t^2}- \Laplace v(x,t)  = d(\mathbf{x},t) \delta(\mathbf{x}-\mathbf{x_r}) \\
    & v(\mathbf{x}, T) = 0                \\
    & v_t(\mathbf{x}, T ) = 0                \\
     \end{array} \right.
\end{equation}
Here $T$ is the final recording time, $d$ is the observed data from the receiver $\mathbf{x_r}$ and $m$ is the assumed background velocity. The adjoint wave equation~\eqref{eq:rtm_adj} is always solved backward in time from $T$ to 0. Therefore it is also referred as backward propagation.

\begin{figure}
\centering
	\subfloat[]{\includegraphics[height=0.2\textwidth]{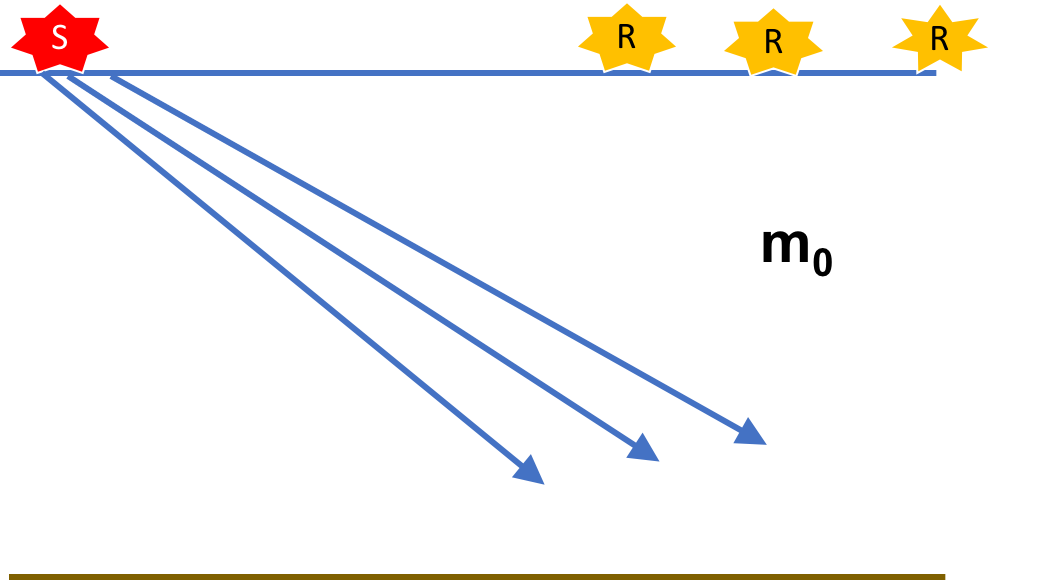}
\label{fig:RTM1} }
	\subfloat[]{\includegraphics[height=0.2\textwidth]{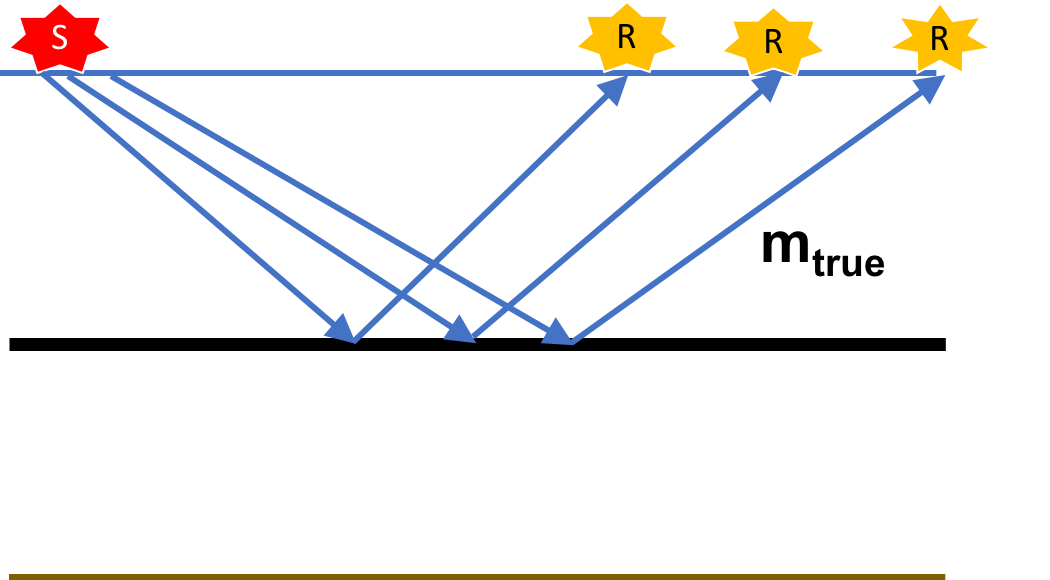}\label{fig:RTM2}}
		\subfloat[]{\includegraphics[height=0.2\textwidth]{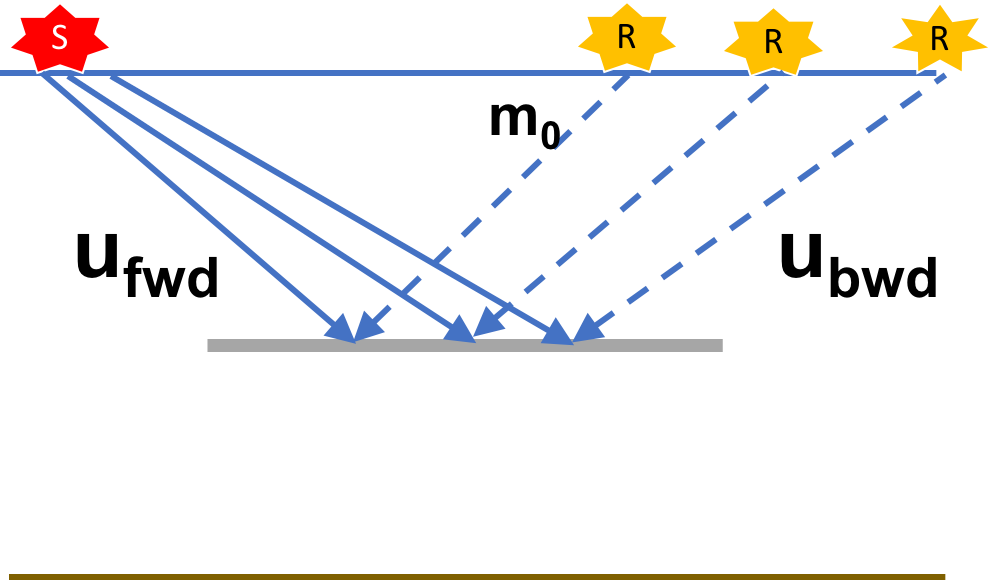}\label{fig:RTM3}}
\caption{RTM: (a)~Synthetic forward wavefield $u_\text{fwd}$, (b)~True forward wavefield and (c)~Reflectors generated as the backward wavefield $u_\text{bwd}$ cross-correlated with $u_\text{fwd}$}\label{fig:RTM3step}
\end{figure}

In classical RTM, the forward modeling typically does not contain reflection information. For example, it can be the paraxial approximation of the wave equation, which does not allow for reflections~\cite{Clayton1977}, or a smooth velocity model with unknown reflecting layers. 
As a summary, the conventional RTM consists three steps as Figure~\ref{fig:RTM3step} shows:
\begin{enumerate}
\item Forward modeling of a wave field with a good velocity model to get $u_\text{fwd}$;
\item Backpropagation of the measured data through the same model to get $u_\text{bwd}$;
\item Cross-correlation the source wavefield $u_\text{fwd}$ and receiver wavefield $u_\text{bwd}$ based on an imaging condition (e.g., Equation~\eqref{eq:IC}) to detect the reflecting interfaces.
\end{enumerate}

RTM uses the entire solution of the wave equations instead of separating the downgoing or upgoing wavefields. Theoretically, RTM produces a more accurate image than ray-based methods since it does not rely on the asymptotic theory or migration using the one-way equation, which typically introduces modeling errors~\cite{Schuster2011}. 
A good background velocity model that contains accurate information about the low-wavenumber components is also crucial for the quality of the image~\cite{Gray2011}. Recent advances in computation power make it possible to compute and store the solution of the wave equation efficiently, which significantly aids RTM to generate high-quality images~\cite{Etgen2009}.

\subsection{Least-squares Reverse Time Migration}
Least-squares reverse time migration (LSRTM) is a new migration method designed to improve the image quality generated by RTM. It is formulated as a linear inverse problem based on the Born approximation which we will describe briefly in this section. The wave equation~\eqref{eq:FWD} defines a nonlinear operator $\bF$ from model domain to data domain that maps $m$ to $u$. The Born approximation is a linearization of this map to the first order so that we can denote it as $L = \frac{\delta \bF}{\delta m}$~\cite{hudson1981use, van1954correlations}.

One can derive the Born approximation as follows~\cite{Demanet2016}. If we denote the model $m(\mathbf{x})$ as the sum of a background model and a small perturbation:
\bq\label{eq:m+dm}
m(\mathbf{x}) = m_0(\mathbf{x}) + \varepsilon m_1(\mathbf{x}),
\eq
the corresponding wavefield $u$ also splits into two parts:
\bq
u(\mathbf{x},t) = u_0(\mathbf{x},t) + u_{sc}(\mathbf{x},t),
\eq
where $u$ satisfies~\eqref{eq:FWD}, and $u_0$ solves the following equation:
\begin{equation}\label{eq:FWD0}
     \left\{
     \begin{array}{rl}
     & m_0(\mathbf{x})\frac{\partial^2 u_0(\mathbf{x},t)}{\partial t^2} - \Laplace u_0(\mathbf{x},t) = s(\mathbf{x},t)\\
    & u_0(\mathbf{x}, 0 ) = 0                \\
    & \frac{\partial u_0}{\partial t}(\mathbf{x}, 0 ) = 0                \\
     \end{array} \right.
\end{equation}

Subtracting \eqref{eq:FWD0} from \eqref{eq:FWD} and using~\eqref{eq:m+dm} , we derive an equation of $u_{sc}$ with zero initial conditions:
\bq \label{eq:scatter}
m_0\frac{\partial^2 u_{sc}(\mathbf{x},t)}{\partial t^2}- \Laplace u_{sc}(\mathbf{x},t)  =  -\varepsilon m_1 \frac{\partial^2 u(\mathbf{x},t)}{\partial t^2}.
\eq

We can write $u_{sc}$ using Green's function $G$:
\bq
u_{sc}(\mathbf{x},t) = -\varepsilon \int_0^t \int_{\R^n} G(\mathbf{x},y;t-s) m_1(y) \frac{\partial^2 u}{\partial t^2} (y,s)dy ds.
\eq
As a result, the original wavefield $u$ has an implicit relation:
\bq
u = u_0 - \varepsilon Gm_1 \frac{\partial^2 u}{\partial t^2} =  \left[I + \varepsilon G m_1  \frac{\partial^2}{\partial t^2}  \right]^{-1} u_0
\eq
The last term can be expanded in terms of Born series,
\begin{eqnarray}
u  &=& u_0 - \varepsilon \int_0^t \int_{\R^n} G(\mathbf{x},y;t-s) m_1(y) \frac{\partial^2 u_0}{\partial t^2} (y,s)dy ds + \bO (\varepsilon ^2)\\
&=& u_0 + \varepsilon u_1 +  \bO (\varepsilon ^2)
\end{eqnarray}

Therefore, we can approximate $u_{sc}$ explicitly by $\varepsilon u_1$ as $- \varepsilon G m_1 \frac{\partial^2 u_0}{\partial t^2}$, which is called the Born approximation. We also derive a linear map from $m_1$ to $u_1$:
\begin{equation}\label{eq:Born}
     \left\{
     \begin{array}{rl}
     & m_0\frac{\partial^2 u_1(\mathbf{x},t)}{\partial t^2} - \Laplace u_1(\mathbf{x},t) = -m_1 \frac{\partial^2 u_0(\mathbf{x},t)}{\partial t^2} \\
    & u_1(\mathbf{x}, 0 ) = 0                \\
    & \frac{\partial u_1}{\partial t}(\mathbf{x}, 0 ) = 0                \\
     \end{array} \right.
\end{equation}
Unlike~\eqref{eq:scatter}, \eqref{eq:Born} is an explicit formulation with $m_0$ as the background velocity and $u_0$ as the background wavefiled which is the solution to~\eqref{eq:FWD0}.

It is convenient to denote the nonlinear forward map~\eqref{eq:FWD} as $\bF: m \mapsto u$. A Taylor expansion of $u = \bF(m)$ in the sense of calculus of variation, gives us:
\bq
u = u_0 + \varepsilon \frac{\delta \bF}{\delta m}[m_0] m_1 + \frac{\varepsilon^2}{2} < \frac{\delta^2 \bF}{\delta m^2}[m_0] m_1, m_1> + \ldots
\eq 
The functional derivative $ \frac{\delta \bF}{\delta m}: m_1 \mapsto u_1 $ is the linear operator~\eqref{eq:Born}, which we hereafter denote as $L$ . The convergence of the Born series and the accuracy of the Born approximation can be proved mathematically~\cite{natterer2004error,newton2013scattering}.

We assume there is an accurate background velocity model $m_0$. The Born modeling operator maps the reflectivity $m_r$ to the scatted wavefield $d_r = \bF(m) - \bF(m_0)$:
\bq
Lm_r = d_r
\eq
Although $L$ is linear, there is no guarantee that it is invertible~\cite{claerbout1992earth}. Instead of computing $L^{-1}$, we seek the reflectivity model by minimizing the least-squares error between observed data $d_r$ and predicted scattering wavefield:
\bq~\label{eq:LSRTM}
J(m_r) = ||Lm_r - d_r||_2^2
\eq
The normal least-squares solution to~\eqref{eq:LSRTM} is $m_r = (L^TL)^{-1}L^Td_r$ where $L^T$  is the adjoint operator, but it is numerically expensive and unstable to invert the term $L^TL$ directly. Instead, the problem is solved in an iterative manner using optimization methods such as conjugate gradient descent (CG). 

Another interesting way of approximating $(L^TL)^{-1}$ is to consider the problem as finding a non-stationary matching filter~\cite{guitton2004amplitude,He2013}. Similar to RTM, we can get an image by doing one step of migration:
\bq~\label{eq:m1}
m_1 = L^Td_r.
\eq
One step of de-migration (Born modeling) based on $m_1$ generates data $d_1$
\bq \label{eq:m2}
d_1 = Lm_1.
\eq
Finally, the re-migration step provides another image $m_2$
\bq  \label{eq:m3}
m_2 = L^Td_1.
\eq
Combining \eqref{eq:m1} to \eqref{eq:m3}, the inverse Hessian operator $(L^TL)^{-1}$ behaves like a matching filter between $m_1$ and $m_2$ which we are able to produce from the observed data. It is also the filter between $m_r$ and $m_1$ as \eqref{eq:m4} and \eqref{eq:m5} show below:
\bq  \label{eq:m4}
m_1 = (L^TL)^{-1} m_2
\eq
\bq \label{eq:m5}
m_r = (L^TL)^{-1} m_1
\eq

Therefore, LSRTM can be seen as a process which first derives a filter to match the re-migration $m_2$ to the initial migration $m_1$ and then applies the filter back to the initial migrated image to give an estimate of the reflectivity. Seeking the reflectivity is equivalent to finding the best filter $K$ by minimizing the misfit $J(K)$ in the image or model domain:
\bq
J(K) = ||m_1 - Km_2 ||_2^2.
\eq
The final reflectivity image $m_r \approx Km_1$. It is a single-iteration method which greatly reduces the computational cost of the iterative methods like CG.

A potentially better way of implementing the filter-based idea is to transform the image into curvelet domain~\cite{candes2006fast} to improve the stability and structural consistency in the matching~\cite{Wang2016}. The formulation of obtaining the Hessian filter in curvelet domain is to minimize a misfit function $J(s)$ where
\bq
J(s) = ||C(m_1) - s C(m_2)||_2^2 + \varepsilon ||s||_2^2,
\eq
where $C$ is the curvelet domain transform operator, $s$ is the matching filter and $\varepsilon$ is the Tikhonov regularization parameter. The final reflectivity image $m_r \approx C^{-1}(|s|C(m_1))$, where $C^{-1}$ is the inverse curvelet transform operator.

In general, least-squares reverse time migration (LSRTM) is still facing challenges. First of all, the image quality highly depends on the accuracy of the background velocity model $m_0$. Even a small error can make the two wavefields meet at a wrong location, which generates a blurred image or an incorrect reflectivity~\cite{luo2014least}. Another drawback is its high computational cost compared with other traditional migration techniques. In practice, LSRTM fits not only the data but also the noise in the data. Consequently, it boosts the high-frequency noise in the image during the iterative inversion~\cite{Dai2017, Zeng2017}.


%
%

\subsection{Inversion}
The process of imaging through modeling the velocity structure is a form of inversion of seismic data~\cite{Treitel2001}, but in this paper, we regard inversion as a process of recovering the quantitative features of the geographical structure, that is, finding $m(\mathbf{x})$ in~\eqref{eq:FWD}. Inversion is often used to build a velocity model iteratively until the synthetic data matches the actual recording~\cite{Mora1989}.

Wave equation traveltime tomography~\cite{luo1991wave} and the ray-based tomography in the earlier section are phase-like inversion methods~\cite{Schuster2011}.  Least-squares inversion is known as linearized waveform inversion~\cite{lailly1984migration,tarantola1984linearized}. The migration method introduced earlier, LSRTM, can also be seen as a linear inverse problem. The background model $m_0$ is not updated after each iteration in least-squares inversion. Similar to the goal of migration, the model to be updated iteratively is the reflectivity distribution instead of the velocity model. One can interpret the process as a series of reverse time migrations, where the data residual is backpropagated into the model instead of the recorded data itself (Figure~\ref{fig:RTM3}).

If the background model $m_0$ is the parameter we invert for, the problem turns into a nonlinear waveform inversion, which is also called full-waveform inversion (FWI). Both the low-wavenumber and high-wavenumber components are updated simultaneously in FWI so that the final image has high resolution and high accuracy~\cite{virieux2009overview}. FWI is the primary focus of the paper. In the following sections, we will further discuss the topic and especially the merit of using optimal transport based ideas to tackle the current limitations. 

\begin{figure}
\includegraphics[scale=0.5]{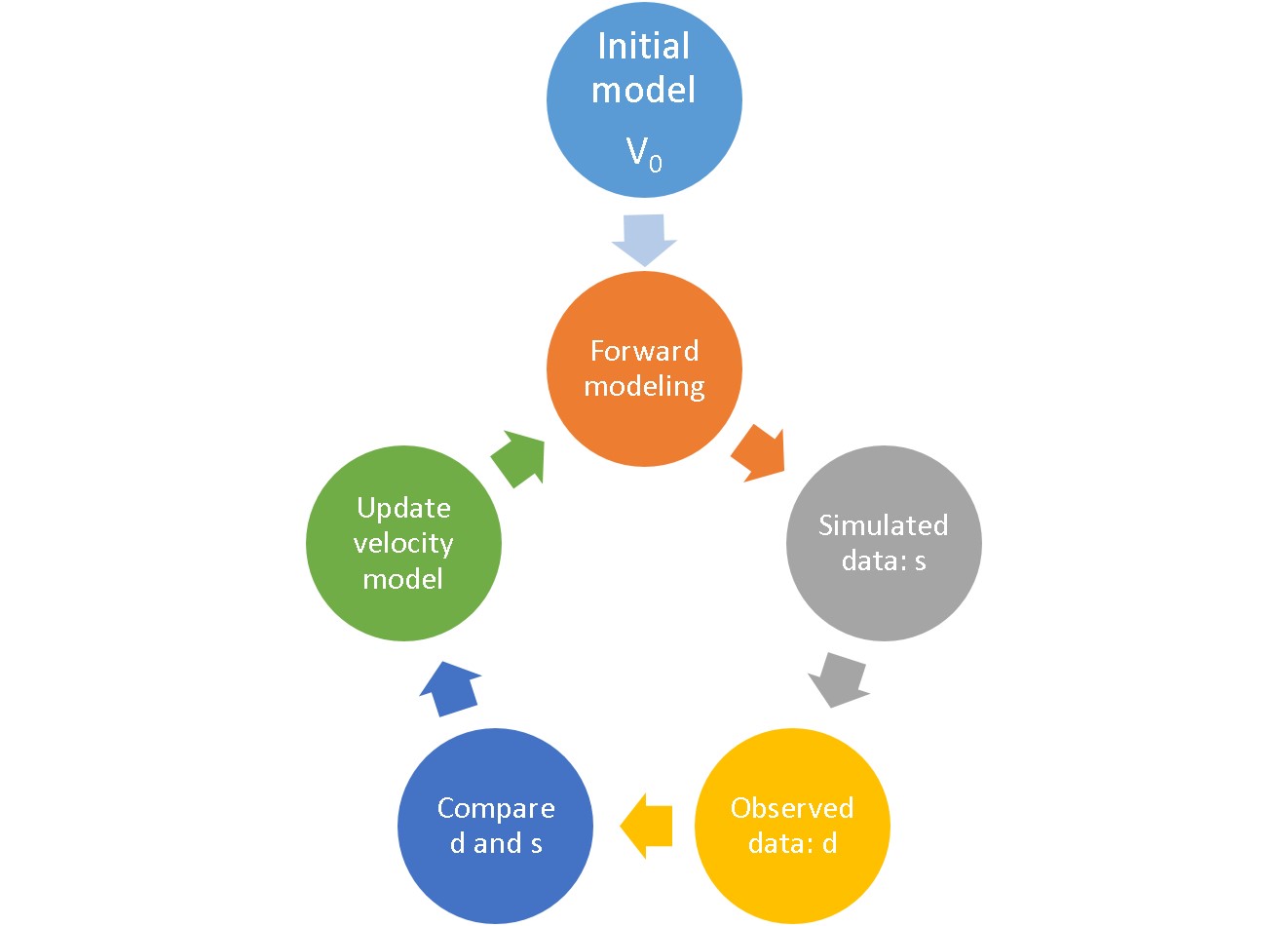}
\caption{The framework of FWI as a PDE-constrained optimization} \label{fig:FWI}
\end{figure}

\section{Full Waveform Inversion}
FWI is a nonlinear inverse technique that
utilizes the entire wavefield information to estimate
the earth properties. 
The notion of FWI was first brought up three decades ago~\cite{lailly1983seismic, tarantola1982generalized} and has been actively studied as the computing power increases.
As we will see, the mathematical formulation of FWI is PDE constrained optimization.
Even inversion for subsurface elastic parameters using FWI
has become increasingly popular in exploration applications~\cite{ brossier2009seismic, Mora1988, virieux2009overview}. 
Currently, FWI can achieve stunning clarity and resolution. Both academia and industry have been actively working on the innovative algorithms and software of FWI. However, this technique is still facing three main challenges.

First, the physics of seismic waves are complex, and we need more accurate forward modeling in inversion going from pure acoustic waves to anisotropic viscoelasticity. Recent developments focus on this multiparameter and multi-mode modeling. FWI strategies for simultaneous and hierarchical velocity and attenuation inversion have been investigated recently~\cite{Qu2017}, but there is a dilemma. The more realistic with more parameters the models of the earth become, the more ill-posed and even non-unique will the inverse problem be.

Second, it is well known that the accuracy of FWI deteriorates from the lack of low frequencies, data noise, and poor starting model. The limitation is mainly due to the ill-posedness of the inverse problem which we treat as a PDE-constrained optimization. FWI is typically performed using local optimization methods in which the subsurface model is described by using a large number of unknowns, and the number of model parameters is determined a priori~\cite{tarantola2005inverse}. These methods typically only use the local gradient of the objective function. As a result, the inversion process is easily trapped in the local minima. Markov chain Monte Carlo (MCMC) based methods~\cite{Sambridge2011}, particle swarm optimization~\cite{Chen2017}, and many other global optimization methods~\cite{Sen2011} can avoid the pitfall theoretically, but they are not cost-efficient to handle practical large-scale inversion currently.

Third, it is relatively inexpensive to update the model through local optimization methods in FWI, but the convergence of the algorithm highly depends on the choice of a starting model. The research directions can be grouped into two main ideas to tackle this problem. One idea is to replace the conventional least-squares norm with other objective functions in optimization for a wider basin of attraction~\cite{engquist2016optimal}. The other idea is to expand the dimensionality of the unknown model by adding non-physical coefficients. The additional coefficients may convexify the problem and fit the data better~\cite{ biondi2012tomographic, huang2017full}. 

The essential elements of FWI framework (Figure~\ref{fig:FWI}) includes forward modeling and the adjoint-state method for gradient calculation.

\subsection{Forward modeling}
Wave-propagation modeling is the most significant step in seismic imaging. The earth is complex with various heterogeneity on many scales, and the real physics is far more complicated than the simple acoustic setting of this paper, 
but the industry standard is still the acoustic model in time or frequency domain.
The current research of FWI covers multiple parameters inversion of seismic waveforms including anisotropic parameters, density, and attenuation factors~\cite{yang2016review} including viscoelastic modeling which is related to fractional Laplacian wave equations~\cite{schiessel1995generalized}. It should be noted that the more parameters in a model, the less well-posed is the inverse problem.

If we exclude the attenuation parameter, the general elastic wave equation is a realistic model. Based on the equation of conservation
of momentum (Newton's law of dynamics) and Hooke's law for stress and strain tensors, we have the following elastic wave equation:
\bq\label{eq:elastic1}
\rho \frac{\partial^2 u_i}{\partial t^2} = f_i + \frac{\partial \sigma_{ij}}{\partial x_j},
\eq
\bq\label{eq:elastic2}
\frac{\partial \sigma_{ij}}{\partial t} = c_{ijkl} \frac{\partial \varepsilon_{ij}}{\partial t}   + \frac{\partial \widetilde{\sigma}_{ij}}{\partial t},
\eq
where $\rho$ is the density, $\mathbf{u}$ is the displacement vector, $\sigma$ is the nine-component stress tensor (i,j = 1,2,3), $\widetilde{\sigma}$ is the internal stress, $\mathbf{f}$ is the outer body force, $\varepsilon$ is the nine-component strain tensor which satisfies $\varepsilon_{ij} = \frac{1}{2} \left( \frac{\partial u_i}{\partial x_j} + \frac{\partial u_j}{\partial x_i} \right)$ and $c_{ijkl}$ is the stiffness tensor containing twenty-one independent components.

One can classify the current numerical methods of complex wave propagation into three categories: direct methods, integral wave equation methods and asymptotic methods~\cite{Viveros2011}. Direct methods include finite-difference method (FDM)~\cite{moczo2007finite}, pseudospectral method~\cite{furumura1998parallel}, finite element method (FEM)~\cite{marfurt1984accuracy}, spectral element method (SEM)~\cite{komatitsch1999introduction}, discontinuous Galerkin method (DG)~\cite{kaser2006arbitrary}, etc. Integral wave equation methods include both boundary element method (BEM)~\cite{bouchon2007boundary} and the indirect boundary element methods (IBEM)~\cite{pointer1998numerical} with a fast multipole method (FMM)~\cite{fujiwara2000fast} for efficiency. Asymptotic methods include geometrical optics, Gaussian beams~\cite{vcerveny1982computation} and frozen Gaussian beams~\cite{lu2011frozen}.

\subsection{Measure of mismatch}
In seismic inversion, the misfit function, i.e. the objective function in the optimization process, is defined as a functional on the data domain. Common misfit functions include cross-correlation traveltime measurements~\cite{luo1991wave,Marquering1999}, amplitude variations~\cite{Dahlen2002} and waveform differences~\cite{tarantola1982generalized}. In both time~\cite{tarantola1987inverse} and frequency domain~\cite{PRATT1990,pratt1990inverse}, the least-squares norm has been the most widely used misfit function. For example, in time domain conventional FWI defines a least-squares waveform misfit as
\begin{equation}
d(f,g) = J(m)=\frac{1}{2}\sum_r\int\abs{f(\mathbf{x_r},t;m)-g(\mathbf{x_r},t)}^2dt,
\end{equation}
where $\mathbf{x_r}$ are receiver locations, $g$ is observed data, and $f$ is simulated data which solves~\eqref{eq:FWD} with model parameter $m$. The time integral is carried out numerically as a sum. This formulation can also be extended to the case with multiple sources.

Real seismic data usually contains noise. As a result, denoising becomes an important step in seismic data processing. The $L^2$ norm is well known to be sensitive to noise~\cite{brossier2010data}. Other norms have been proposed to mitigate this problem. For example, the $L^1$ norm~\cite{Crase1990,tarantola1987inverse}, the Huber criterion~\cite{Guitton2003,ha2009waveform} and the hybrid $L^1/L^2$ criterion~\cite{Bube1997} all demonstrated improved robustness to noise compared with conventional $L^2$ norm. 

All the misfit functions above are point-by-point based objective functions which means they only accumulate the differences in amplitude at each fixed time grid point. There are global misfit functions that compare the simulated  and measured signals not just pointwise. The Wasserstein metric is one such metric which we will discuss later. It is very robust with respect to noise

\begin{figure}
	\subfloat[Two signals]{\includegraphics[width=1.0\textwidth]{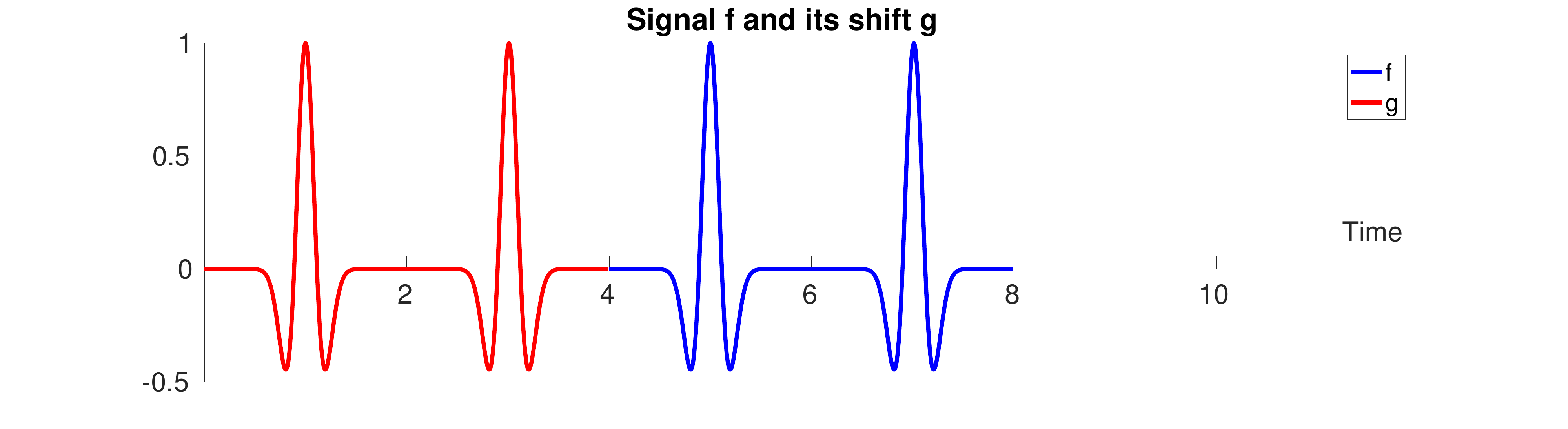}\label{fig:2_ricker_signal}}\\
  	\subfloat[$L^2$ sensitivity curve]
{\includegraphics[width=0.5\textwidth]{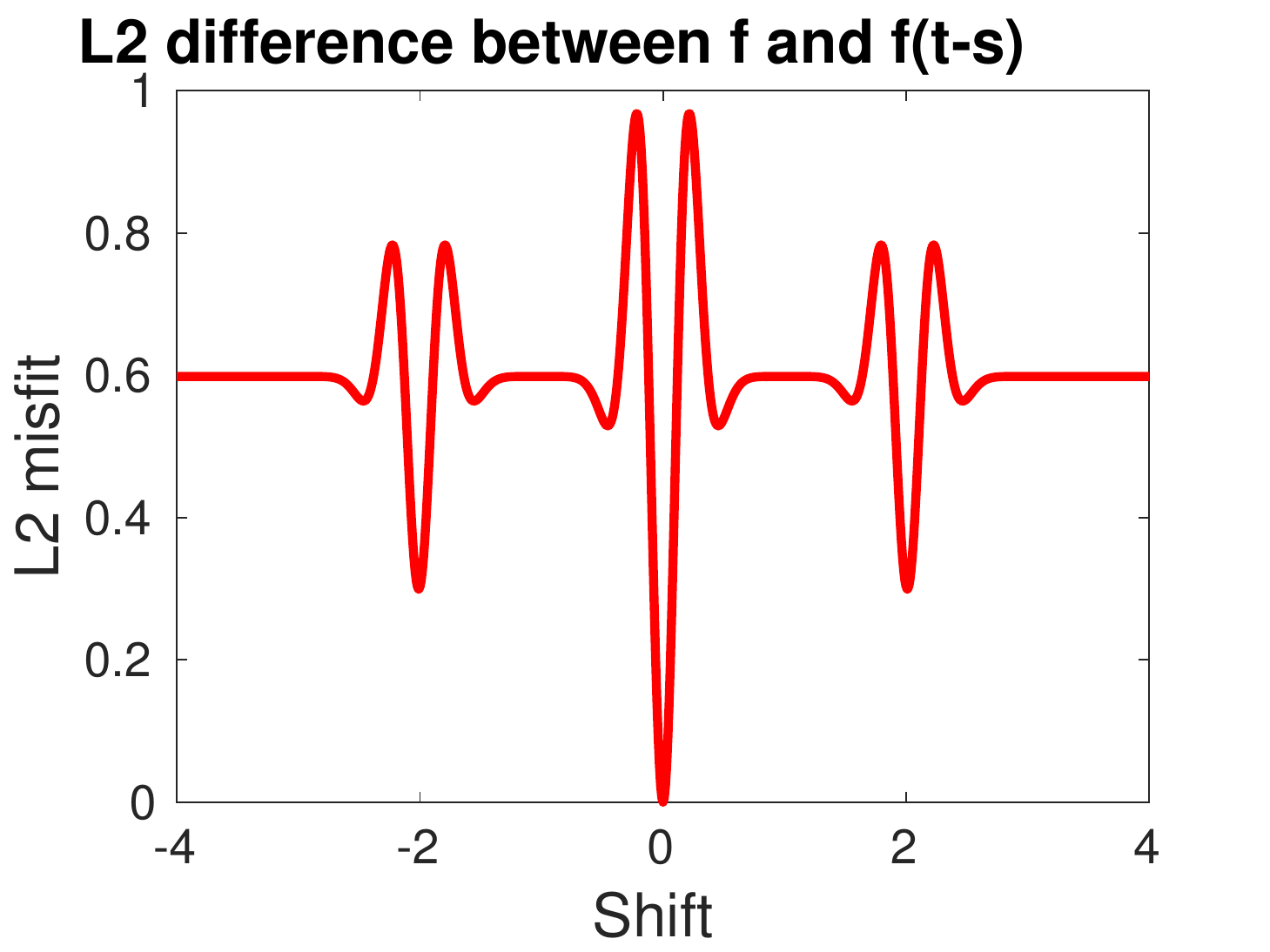}\label{fig:2_ricker_L2}}
  	\subfloat[$W_2$ sensitivity curve]
{\includegraphics[width=0.5\textwidth]{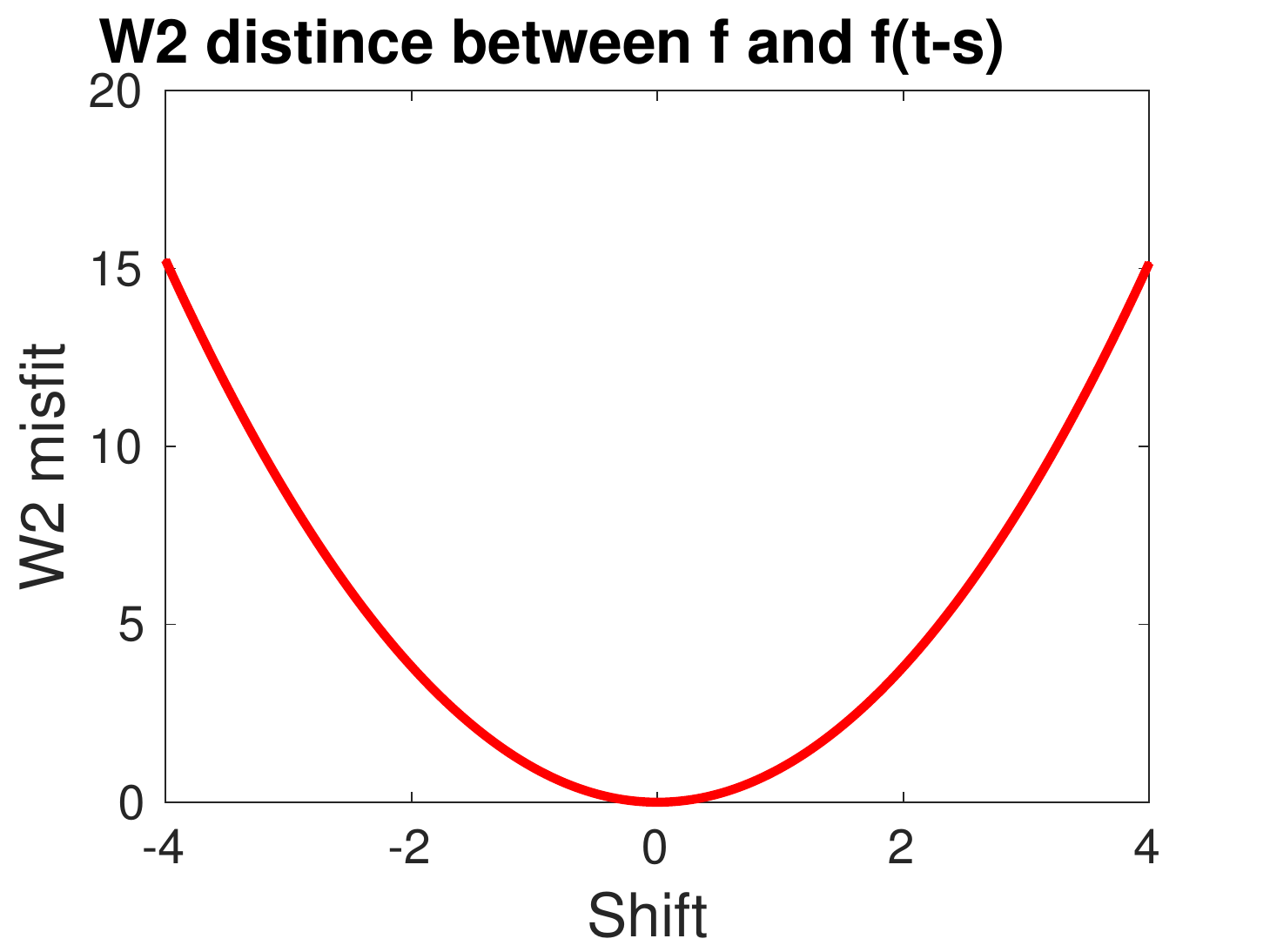}\label{fig:2_ricker_W2}}

\caption{(a)~A signal consisting two Ricker wavelets (blue) and its shift (red)~(b)~$L^2$ norm of the difference between $f$ and $f(t-s)$ in terms of shift $s$~(c)~$W_2$ norm between $f$ and $f(t-s)$ in terms of shift $s$}
\end{figure}

The oscillatory and periodic nature of waveforms lead to another main challenge in FWI: the cycle-skipping issue when implementing FWI as a local inversion scheme. 
If the true data and the initial synthetic data are more than half wavelength ($>\frac{\lambda}{2}$) away from each other, the first gradient can go in the wrong direction regarding the phase mismatch, but can nonetheless reduce the data misfit in the fastest manner~\cite{Beydoun1988}. Mathematically, it is related to the highly nonconvex and highly nonlinear nature of the inverse problem and results in finding only a local minima.
Figure~\ref{fig:2_ricker_signal} displays two signals, each of which contains two Ricker wavelets and $f$ is simply a shift of $g$. The $L^2$ norm between $f$ and $g$ is plotted in Figure~\ref{fig:2_ricker_L2} as a function of the shift $s$. We observe many local minima and maxima in this simple two-event setting which again demonstrated the difficulty of the, so called,  cycle-skipping issues~\cite{yang2017analysis}. 

\begin{figure}	
	\centering
  	\subfloat[Misfit for $L^2$ norm]{\includegraphics[width=0.45\textwidth]{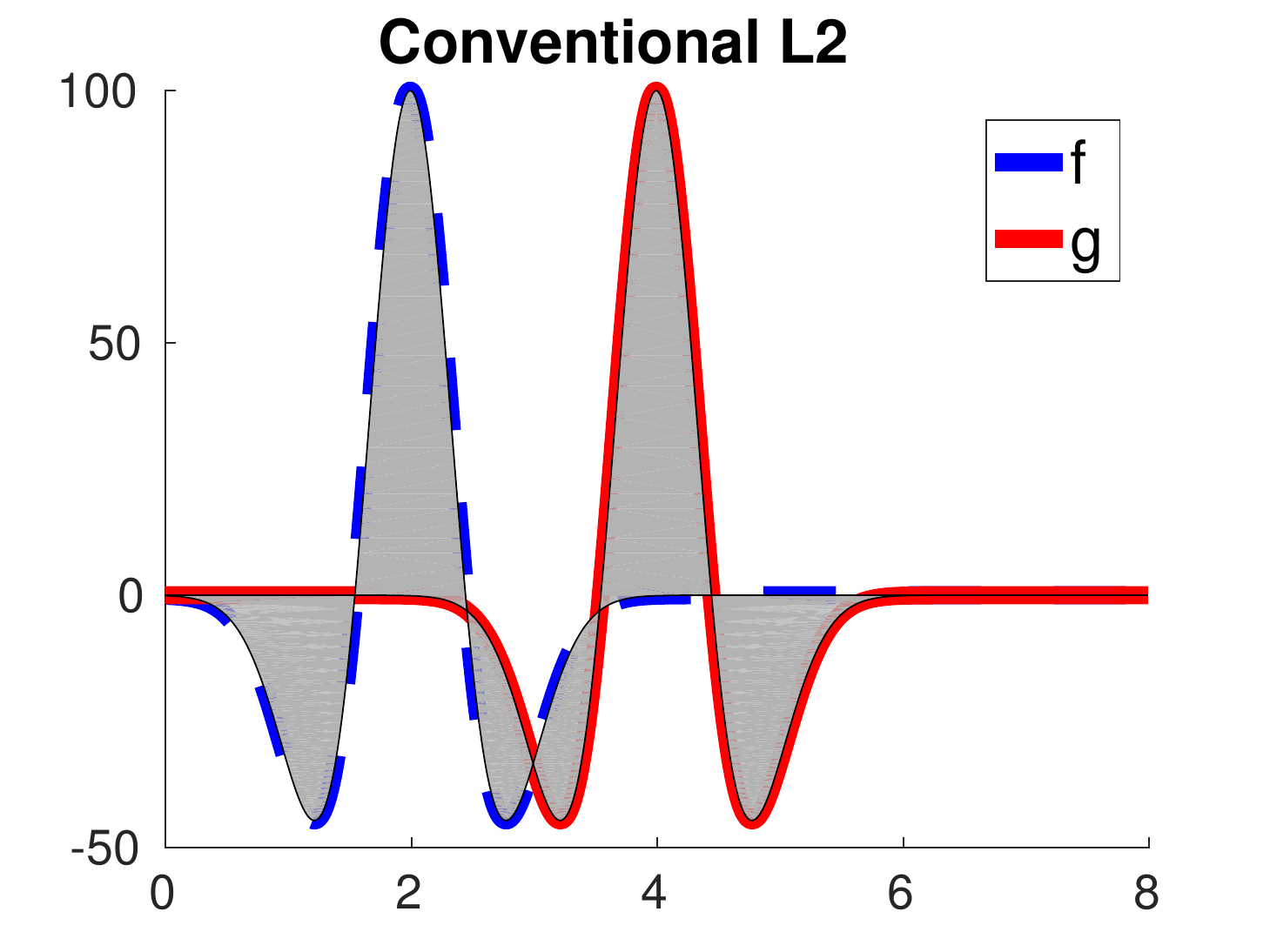}\label{fig:L2(fg)}}
  	\subfloat[Misfit for Integral $L^2$ method]{\includegraphics[width=0.45\textwidth]{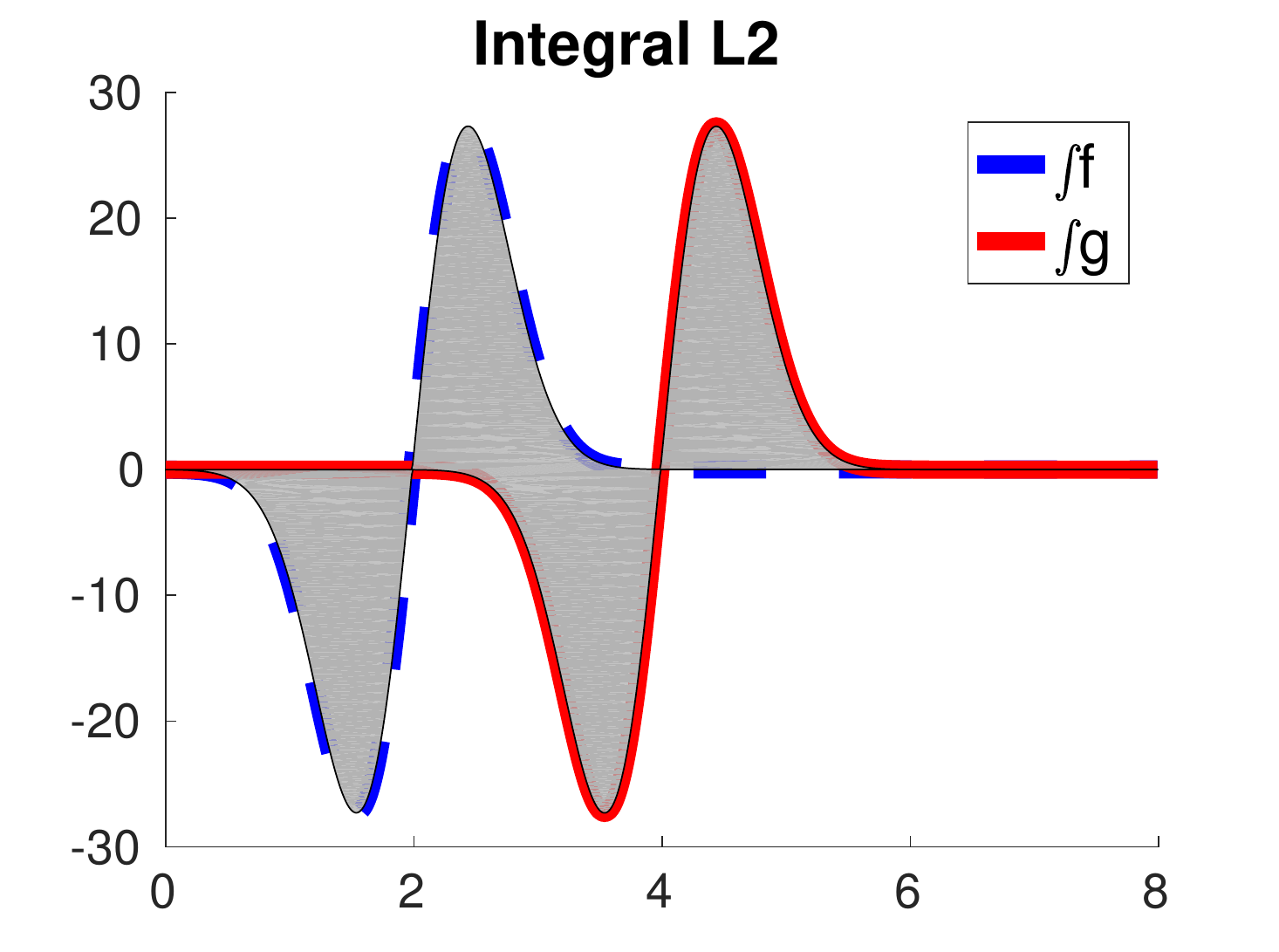}\label{fig:I(fg)}} 
	\caption{The shaded areas represent the mismatch each misfit function considers. (a) $L^2$: $\int (f-g)^2 dt$. (b) Integral wavefields method: $\int (\int f-\int g)^2 dt$.~\cite{yangletter}}
	\label{fig:2L2}
\end{figure}

The lower frequency components have a wider basin of attraction with the least-squares norm being the misfit function. Several hierarchical methods that invert from low frequencies to higher frequencies have been proposed in the literature to mitigate the cycle-skipping of the inverse problem~\cite{bunks1995multiscale, kolb1986pre, pratt1990inverse,sirgue2004efficient,weglein2003inverse}. Several other methods instead compare the integrated waveforms~\cite{huang2014two, liu2012normalized} (Figure~\ref{fig:2L2}) and the waveform envelops~\cite{Bozdag2011, luo2015seismic}. They share a similar idea with the hierarchical methods of taking advantage of the lower frequency components in the data.

A recently introduced class of misfit functions is based on optimal transport~\cite{chen2017quadratic,EFWass, engquist2016optimal,W1_2D,W1_3D,yang2017analysis,yangletter, yang2017application}. As a useful tool from the theory of optimal transport, the Wasserstein metric computes the minimal cost of rearranging one distribution into another. The optimal transport based methods compare the observed and simulated data globally and thus include phase information. We will discuss these measures in section 4 and 5.

Other misfit functions with the idea of non-local comparison proposed in the literature include filter based misfit functions \cite{AWI,zhu2016building} as well as inversion using, so called, dynamic time warping~\cite{ma2013wave} and the registration map~\cite{Baek}. The differential semblance optimization~\cite{symes1991velocity} exploits both phase and amplitude information of the reflections. Tomographic full waveform inversion \cite{biondi2012tomographic} has some global convergence characteristics of wave-equation migration velocity analysis. 
In the filter based methods~\cite{AWI,zhu2016building}, a filter is designed to minimize the $L^2$ difference between filtered simulated data and the observed data. The misfit is then a measure of how much the filter deviates from the identity. As we will see in the optimal transport based technique, this is done in one step where the optimal map directly determines the mapping of the simulated data. The optimal transport map is general and does not need to have the form of a convolution filter as in the filter based methods. 

\subsection{Adjoint-state method}
Large-scale realistic 3D inversion is possible today. The advances in numerical methods and computational power allow for solving the 3D wave equations and compute the Fr\'{e}chet derivative with respect to model parameters, which are needed in the optimization.  In the adjoint-state method, one only needs to solve two wave equations numerically, the forward propagation and the backward adjoint wavefield propagation. Different misfit functions typically only affect the source term in the adjoint wave equation~\cite{Plessix,tarantola2005inverse}. 

Let us consider the misfit function $J(m)$ for computing the difference between predicted data $f$ and observed data $g$ where $m$ is the model parameter, $F(m)$ is the forward modeling operator, $u(\mathbf{x},t)$ is the wavefield and $s(\mathbf{x},t)$ is the source. The predicted data $f$ is the partial Cauchy boundary data of $u$ which can be written as $f = Ru$ where $R$ is a restriction operator only at the receiver locations. The wave equation~\eqref{eq:FWD} can be denoted as
\bq ~\label{eq:adj_fwd}
F(m) u = s.
\eq
Taking first derivative regarding model $m$ on both sides gives us:
\bq
\frac{\partial F}{\partial m} u + F \frac{\partial u}{\partial m} = 0.
\eq
Therefore,
\bq~\label{eq:adj_grad1}
\frac{\partial f}{\partial m}  = -RF^{-1} \frac{\partial F}{\partial m} u.
\eq

By the chain rule, the gradient of misfit function $J$ with respect to $m$ is
\bq~\label{eq:adj_grad0}
\frac{\partial J}{\partial m} = \left(\frac{\partial f}{\partial m} \right)^T  \frac{\partial J}{\partial f}
\eq
We can derive the following equation by plugging~\eqref{eq:adj_grad1} into ~\eqref{eq:adj_grad0} :
\bq~\label{eq:adj_grad2}
\frac{\partial J}{\partial m}  =-u^T \left(\frac{\partial F}{\partial m}\right)^T F^{-T}R^T  \frac{\partial J}{\partial f}
\eq

Equation~\eqref{eq:adj_grad2} is the adjoint-state method. The term $F^{-T}R^T  \frac{\partial J}{\partial f}$ denotes the backward wavefield $v$ generated by the adjoint wave equation whose source is the data residual $R^T  \frac{\partial J}{\partial f}$. The gradient is similar to the usual imaging condition~\eqref{eq:IC}:
\bq~\label{eq:adj_grad3}
\frac{\partial J}{\partial m}  =- \int_0^T \frac{\partial^2 u(\mathbf{x},t)}{\partial t^2} v(\mathbf{x},t)dt,
\eq
where $v$ is the solution to the adjoint wave equation:
\begin{equation} \label{eq:FWI_adj}
     \left\{
     \begin{array}{rl}
     & m\frac{\partial^2 v(\mathbf{x},t)}{\partial t^2}- \Laplace v(\mathbf{x},t)  = R^T\frac{\partial J}{\partial f}\\
    & v(\mathbf{x}, T) = 0                \\
    & v_t(\mathbf{x}, T ) = 0                \\
     \end{array} \right.
\end{equation}
Therefore $F^T$ can be seen as the backward modeling operator which is similar to the adjoint wave equation~\eqref{eq:rtm_adj} but with a different source term.

There are many other equivalent ways to formulate the adjoint-state method. One can refer to~\cite{Demanet2016,Plessix} for more details.

In FWI, our aim is to find the model parameter $m^{\ast}$ that minimizes the objective function, i.e. \(m^{\ast} = \argmin J(m) \). For this PDE-constrained optimization, one can use the Fr\'{e}chet derivative in a gradient-based iterative scheme to update the model $m$, such as steepest descent, conjugate gradient descent (CG), L-BFGS, Gauss-Newton method, etc. One can also derive the second-order adjoint equation for the Hessian matrix and use the full Newton's method in each iteration, but it is not practical regarding memory and current computing power. It is one of the current research interests to analyze and approximate the Hessian matrix in optimization~\cite{Virieux2017}.

\section{Optimal Transport for FWI}
Optimal transport has become a well-developed topic in mathematics since it was first brought up by Monge~\cite{Monge} in 1781. 
Due to its ability to incorporate both intensity and spatial information, optimal transport based metrics for modeling and signal processing have recently been adopted in a variety of applications including image retrieval, cancer detection, and machine learning \cite{kolouri2016transport}. In computer science, the metric is often called the ``Earth Mover's Distance'' (EMD).

The idea of using optimal transport for seismic inversion was first proposed in~\cite{EFWass}. The Wasserstein metric is a concept based on optimal transportation~\cite{Villani}. Here, we transform our datasets of seismic signals into density functions of two probability distributions. Next, we find the optimal map between these two datasets and compute the corresponding transport cost as the misfit function in FWI.  In this paper, we will focus on the quadratic cost function. The corresponding misfit is the quadratic Wasserstein metric ($W_2$). As Figure~\ref{fig:2_ricker_W2} shows, the convexity of $W_2$ is much better than the $L^2$ norm when comparing oscillatory seismic data with respect to shift.

Following the idea that changes in velocity cause a shift or ``transport''  in the arrival time, \cite{engquist2016optimal} demonstrated the advantageous mathematical properties of the quadratic Wasserstein metric ($W_2$) and provided rigorous proofs that laid a solid theoretical foundation for this new misfit function. We can apply $W_2$ as misfit function in two different ways: trace-by-trace comparison which is related to 1D optimal transport in the time dimension, and the entire dataset comparison in multiple dimensions. We will see that solving the \MA equation in each iteration of FWI is a useful technique~\cite{yang2017application} for calculating the Wasserstein distance. An analysis of the 1D optimal transport approach and the conventional misfit functions such as $L^2$ norm and integral $L^2$ norm illustrated the intrinsic advantages of this transport idea~\cite{yangletter}. 

\subsection{Wasserstein metric}

Let X and Y be two metric spaces with nonnegative Borel measures $\mu$ and $\nu$ respectively. Assume X and Y have equal total measure:
\bq
\int_X d\mu = \int_Y d\nu
\eq
Without loss of generality, we will hereafter assume the total measure to be one, i.e., $\mu$ and $\nu$ are probability measures.

\begin{definition}[Mass-preserving map]
A transport map $T: X \rightarrow Y$ is mass-preserving if for any measurable
set $B \in  Y$ ,
\bq
\mu (T^{-1}(B)) = \nu(B)
\eq
If this condition is satisfied, $\nu$ is said to be the push-forward of $\mu$ by $T$, and we write $\nu = T_\# \mu $
\end{definition}
 
In another word, given two nonnegative densities $f = d\mu$ and $g=d\nu$, we are interested in the mass-preserving map $T$ such that $f = g \circ T$. The transport cost function $c(x,y)$ maps pairs $(x,y) \in X\times Y$ to $\mathbb{R}\cup \{+\infty\}$, which denotes the cost of transporting one unit mass from location $x$ to $y$. The most common choices of $c(x,y)$ include $|x-y|$ and $|x-y|^2$, which denote the Euclidean norms for vectors $x$ and $y$ hereafter. Once we find a mass-preserving map $T$, the cost corresponding to $T$ is 
\[
I(T,f,g,c) = \int\limits_Xc(x,T(x))f(x)\,dx. 
\]

While there are many maps $T$ that can perform the relocation, we are interested in finding the optimal map that minimizes the total cost
\[
I(f,g,c) = \inf\limits_{T\in\M}\int\limits_Xc(x,T(x))f(x)\,dx,
\]
where $\M$ is the set of all maps that rearrange $f$ into $g$.

Thus we have informally defined the optimal transport problem, the optimal map as well as the optimal cost, which is also called the Wasserstein distance:
\begin{definition}[The Wasserstein distance]
  We denote by $\mathscr{P}_p(X)$ the set of probability measures with finite moments of order $p$. For all $p \in [1, \infty)$,   
\bq~\label{eq:static}
W_p(\mu,\nu)=\left( \inf _{T_{\mu,\nu}\in \mathcal{M}}\int_{\mathbb{R}^n}\left|x-T_{\mu,\nu}(x)\right|^p d\mu(x)\right) ^{\frac{1}{p}},\quad \mu, \nu \in \mathscr{P}_p(X).
\eq
$\mathcal{M}$ is the set of all maps that rearrange the distribution $\mu$ into $\nu$.
\end{definition}

\subsection{1D problem}
In~\cite{yang2017application}, we proposed two ways of using $W_2$ in FWI were proposed. One can either compute the misfit globally by solving a 2D or 3D optimal transport problem or compare data trace-by-trace with the 1D explicit formula, see Theorem 1 below. For the 1D approach, the corresponding misfit function in FWI becomes
\bq \label{eqn:Wp1D}
J_1(m) = \sum\limits_{r=1}^R W_2^2(f(\mathbf{x_r},t;m),g(\mathbf{x_r},t)), 
\eq
where $R$ is the total number of time history traces, $g$ is the observed data, $f$ is the simulated data, $\mathbf{x_r}$ are the receiver locations, and $m$ is the model parameter. 
Mathematically it is $W_2$ metric in the time domain and $L^2$ norm in the spatial domain.

For $f$ and $g$ in one dimension, it is possible to exactly solve the optimal transportation problem~\cite{Villani} in terms of the cumulative distribution functions
\bq \label{eq:F&G}
F(x) = \int_{-\infty}^x f(t)\,dt, \quad G(y) = \int_{-\infty}^y g(t)\,dt.
\eq

In fact, the optimal map is just the unique monotone rearrangement of the density $f$ into $g$. In order to compute the Wasserstein metric ($W_p$), we need the cumulative distribution functions $F$ and $G$ and their inverses $F^{-1}$ and $G^{-1}$ as the following theorem states:

\begin{theorem}[Optimal transportation on $\R$]\label{OT1D}
Let $0 < f, g < \infty$ be two probability density functions, each supported on a connected subset of $\R$.  Then the optimal map from $f$ to $g$ is $T = G^{-1}\circ F$.
\end{theorem}


From the theorem above, we derive another formulation for the 1D quadratic Wasserstein metric:
\bq\label{myOT1D}
\begin{aligned}
W_2^2(f,g) & = \int_0^1 |F^{-1} - G^{-1}|^2  dy  \\
& = \int_X |x-G^{-1}(F(x))|^2 f(x)dx.
\end{aligned}
\eq

The corresponding Fr\'{e}chet derive which is also the adjoint source term in the backward propagation is:
\bq~\label{eq:1D_ADS_C}
\begin{split}
\frac{\partial W_2^2(f,g)}{\partial f} = & \left( \int_t^{T_0}-2(s-G^{-1}(F(s))\frac{dG^{-1}(y)}{dy}\biggr\rvert_{y=F(s)}  f(s) ds \right) dt \\ & +
 |t-G^{-1}(F(t))|^2 dt.
\end{split}
\eq

This adjoint source term in the discrete 1D setting can be computed  as
\bq\label{eq:1D_ADS_D}
\begin{split}
 \left[U\ \diag \left(\frac{-2 f(t) dt}{g(G^{-1}\circ F (t))}  \right)   \right] (t-G^{-1} \circ F(t)) dt + |t-G^{-1}\circ F(t)|^2  dt,
\end{split}
\eq
where $U$ is the upper triangular matrix whose non-zero components are 1.

\subsection{\MA equation}
This fully nonlinear partial differential equation plays an important role in computing the Wasserstein metric.
\subsubsection{Introduction}
In the previous section, we introduced the 1D optimal transport technique of comparing seismic data trace by trace and the explicit solution formula. Another option is a general optimal transport problem in all dimensions. In the global case we compare the full datasets and consider the whole synthetic data $f$ and observed data $g$ as objects with the general quadratic Wasserstein metric ($W_2$):
\bq \label{eqn:W22D}
J_2(m) = W_2^2(f(\mathbf{x_r},t;m),g(\mathbf{x_r},t)).
\eq

The simple exact formula for 1D optimal transportation does not extend to optimal transportation in higher dimensions.  Nevertheless, it can be computed by relying on two important properties of the optimal mapping~$T(x)$: conservation of mass and cyclical monotonicity.  From the definition of the problem, $T(x)$ maps $f$ into $g$.  If $T$ is a sufficiently smooth map and $\det(\nabla T(x) ) \neq 0$, the change of variables formula formally leads to the requirement
\bq\label{eq:massConserved}
f(x) = g(T(x))\det(\nabla T(x)).
\eq 

The optimal map takes on additional structure in the special case of the cost function (i.e., $c(x,y) = |x-y|^2$): it is cyclically monotone~\cite{Brenier,KnottSmith}.
\begin{definition}[Cyclical monotonicity]
\label{cyclical}
We say that $T:X\to Y$ is cyclically monotone if for any $m\in\mathbb{N}^+$, \(x_i\in X, \,1\leq i \leq m\),
\bq\label{eq:cyclical}
\sum_{i=1}^{m}|x_i-T(x_i)|^2 \leq  \sum_{i=1}^{m}|x_i-T(x_{i-1})|^2
\eq 
or equivalently 
\bq
\sum_{i=1}^{m}\langle T(x_i),x_i-x_{i-1}\rangle \geq 0  
\eq 
where $x_0\equiv x_m$.
\newline
\end{definition}

Additionally, a cyclically monotone mapping is formally equivalent to the gradient of a convex function~\cite{Brenier,KnottSmith}.  Making the substitution $T(x) = \nabla u(x)$ into the constraint~\eqref{eq:massConserved} leads to the \MA equation
\bq~\label{eq:MAA}
\det (D^2 u(x)) = \frac{f(x)}{g(\nabla u(x))}, \quad u \text{ is convex}.
\eq
In order to compute the misfit between distributions $f$ and $g$, we first compute the optimal map $T(x) = \nabla u(x)$ via the solution of this \MA equation coupled to the non-homogeneous Neumann boundary condition 
\bq\label{eq:BC}
\nabla u(x) \cdot \nu = x\cdot \nu, \,\, x \in \partial X.
\eq
The squared Wasserstein metric is then given by
\bq\label{eq:WassMA}
W_2^2(f,g) = \int_X f(x)\abs{x-\nabla u(x)}^2\,dx.
\eq

For the general \MA equation, the uniqueness of the optimal map is not guaranteed. One need to discuss it in the context of a particular cost function and certain hypothesis. For example, the cyclical monotonicity is the key element in the proof of the following Brenier's theorem~\cite{Brenier, DePhilippis2013} which gives an elegant result about the uniqueness of optimal transport map for the quadratic cost $|x-y|^2$:
\begin{theorem}[Brenier's theorem]
Let $\mu$ and $\nu$ be two compactly supported probability measures on $\R^n$. If $\mu$ is absolutely continuous with respect to the Lebesgue measure, then
\begin{enumerate}
\item There is a unique optimal map $T$ for the cost function $c(x,y) = |x-y|^2$. \item There is a convex function $u: \R^n \rightarrow \R$ such that the optimal map $T$ is given by $T(x) = \nabla u(x)$ for $\mu$-a.e. x.
\end{enumerate}
Furthermore, if $\mu(dx) = f(x)dx$, $\nu(dy) = g(y)dy$, then $T$ is differential $\mu$-a.e. and 
\bq
\det (\nabla T(x)) = \frac{f(x)}{g(T(x))}.
\eq
\end{theorem}

We are here considering the connection between the \MA equation and optimal transport where the transport map is geometric in nature.
The \MA equation is of course also known for many other connections to geometry and mathematical physics. Let us mention a few examples. It arises naturally in many problems such as affine geometry~\cite{Cheng1986}, Riemannian geometry~\cite{aubin2013some}, isometric embedding~\cite{Han2006}, reflector shape design~\cite{Wang1996}, etc. In the last century, treatments about this equation mostly came from the geometric problems above~\cite{Caffarelli1985, Cheng1977,Krylov1995, Minkowski1989, trudinger2008monge}.  If we consider the following general \MA equation:
\bq~\label{eq:generalMA}
\det (D^2 u(x)) = f(x, u , Du),
\eq
when $f = K(x) (1 + |Du|^2)^{(n+2)/2}$, the equation becomes the prescribed Gaussian curvature equation~\cite{DePhilippis2013}. In affine geometry, an affine sphere in the graph satisfies the \MA equation~\eqref{eq:generalMA}. The affine maximal surface satisfies a fourth-order equation which is related to the general \MA equation:
\bq
\sum_{i,j = i}^n U^{ij}\partial_{x_i}\partial_{x_j}\left[\det (D^2 u) \right]^{-\frac{n+1}{n+2}} = 0, 
\eq 
where $U^{ij}$ is the cofactor matrix of $D^2 u$~\cite{Trudinger2008}.

\subsubsection{Weak solutions}
Although the \MA equation is a second-order PDE, there is no guarantee that the classical $C^2$ solution always exists. For the generalized \MA equation~\eqref{eq:generalMA} with homogeneous Dirichlet boundary condition $u=0$ on $\partial \Omega$, it is well-known that there exists a classical convex solution $u \in C^2(\Omega) \cup C(\overline{\Omega})$, when $f$ is strictly positive and sufficiently smooth~\cite{Caffarelli1989, Caffarelli1990, Caffarelli1991}. When the assumptions no longer hold, we solve for two types of weak solutions instead: the Aleksandrov solution and the viscosity solution. One can refer to~\cite{gutierrez2016monge} for more details and proofs of the following definitions and theorems.

Let $\Omega$ be the open subset of $\R^d$ and  $u: \Omega \rightarrow \R$. We denote $\mathcal{P}(\R^d)$ as the set of all subsets of $\R^d$.
\begin{definition}
The normal mapping of $u$, or the subdifferential of $u$, is the set-valued mapping $\partial u: \Omega \rightarrow \mathcal{P}(\R^d)$ defined by 
\bq
\partial u(x_0) = \{ p: u(x) \geq u(x_0)+ p\cdot (x - x_0), \text{\   for all } x \in \Omega \}   
\eq
Given $V\in \Omega$, $\partial u(V) = \cup_{x\in V} \partial u(x)$.
\end{definition}

\begin{theorem}[\MA measure]
If $\Omega$ is open and $u \in C(\Omega)$, then the class
\[ \mathcal{S}  = \{ V\subset \Omega:  \partial u(V)  \text{is Lebesgue measurable} \} \]
is a Borel $\sigma$-algebra. The set function $Mu: \mathcal{S} \rightarrow \overline{\R}$ defined by 
\bq~\label{eq:MA_measure}
Mu(V) = |\partial u(V) |
\eq
is a measure, finite on compact sets, called the \MA measure associated with the function $u$.
\end{theorem}
This is a measure generated by the the \MA operator, which naturally defines the notion of the Aleksandrov solution.
\begin{definition}[Aleksandrov solution]
Let $\nu$ be a Borel measure defined on $\Omega$ which is an open and convex subset of $\R^n$. The convex function $u$ is a weak solution, in the sense of Aleksandrov, to the \MA equation
\bq
\det D^2u = \nu \text{\quad in $\Omega$} 
\eq
if the associated \MA measure $Mu$ defined in~\eqref{eq:MA_measure}  is equal to $\nu$.
\end{definition}

Next we state one existence and uniqueness result for the Aleksandrov solution~\cite{awanou2014discrete}.
\begin{theorem}[Existence and uniqueness of the Aleksandrov solution]
Consider the following Dirichlet problem of the \MA equation
\begin{eqnarray} \label{eq:alek}
\det D^2u &= &\nu \text{\quad in $\Omega$} \\
u &=& g \text{\quad on $\partial \Omega$} \nonumber,
\end{eqnarray}
on a convex bounded domain $\Omega \in \R^d$ with boundary $\partial \Omega$. Assume that $\nu$ is a finite Borel measure  and $g \in C(\partial \Omega)$ which can be extended to a convex function $\tilde{g} \in C(\bar{\Omega})$. Then the \MA equation~\eqref{eq:alek} has a unique convex Aleksandrov solution in $C(\overline{\Omega})$.
\end{theorem}

Aleksandrov's generalized solution corresponds to the curvature measure in the theory of convex bodies~\cite{Trudinger2008}.  A finite difference scheme for computing Aleksandrov measure induced by $D^2u$ in 2D was conducted in~\cite{oliker1989numerical} with the solution $u$ comes as a byproduct~\cite{feng2011analysis}.

Another notion of weak solution is the viscosity solution which occurs naturally if $f$ is continuous in \eqref{eq:generalMA}.
\begin{definition} [Viscosity solution]
Let $u\in C(\Omega)$ be a convex function and $f\in C(\Omega)$, $f\geq 0$. The function $u$ is a viscosity subsolution (supersolution) of \eqref{eq:generalMA} in $\Omega$ if whenever convex function $\phi \in C^2(\Omega)$ and $x_0 \in \Omega$ are such that $(u-\phi)(x) \leq (\geq )(u-\phi)(x_0) $ for all $x$ in the neighborhood of $x_0$, then we must have
\[
\det (D^2\phi(x_0)) \leq (\geq) f(x_0).
\]
The function $u$ is a viscosity solution if it is both a viscosity subsolution and supersolution.
\end{definition}

We can relate these two notions of weak solution in the following proposition:
\begin{proposition}
If $u$ is a Aleksandrov (generalized) solution of~\eqref{eq:generalMA} with $f$ continuous, then $u$ is also a viscosity solution.
\end{proposition}

\subsection{Numerical optimal transport in higher dimensions}
In this section, we will summarize some of the current numerical methods for solving the optimal transport problems in higher dimensions. These methods are based on the equivalent or relaxed formulations of the original Monge's problem. In the end, we will introduce a monotone finite difference \MA solver which is proved to converge to the viscosity solution to~\eqref{eq:MAA}~\cite{barles1991convergence,FroeseTransport}. 

\subsubsection{General methods}
Optimal transport is a well-studied subject in mathematics while the computation techniques are comparatively underdeveloped. 
We will focus on analysis based methods. There are combinatorial techniques that typically are computationally costly in higher dimensions, for example, the Hungarian algorithm~\cite{kuhn1955hungarian}.
 
The definition~\eqref{eq:static} is the original static formulation of the optimal transport problem with a quadratic cost. It is an infinite dimensional optimization problem if we search for $T$ directly. The non-symmetric nature of Monge's problem also generated difficulty because the map is unnecessarily bijective~\cite{Levy2017}. 

In the 40's, Kantorovich relaxed the constraints and formulated the dual problem~\cite{kantorovich1960mathematical}. Instead of searching for a map $T$, the transference plan $\gamma$ is considered, which is also a measure supported by the product space $X\times Y$. The Kantorovich problem is the following:
\bq
\inf_{\gamma} \bigg\{ \int_{X \times Y} c(x,y) d\gamma\ |\ \gamma \geq 0\ \text{and}\ \gamma \in \Pi(\mu, \nu) \bigg\} ,
\eq
where $\Pi (\mu, \nu) =\{ \gamma \in \mathcal{P}(X\times Y)\ |\ (P_X)\sharp \gamma = \mu, (P_Y)\sharp \gamma = \nu  \}$. Here $(P_X)$ and $(P_Y)$ denote the two projections, and $(P_X)\sharp \gamma$ and $(P_Y)\sharp \gamma $ are two measures obtained by pushing forward $\gamma$ with these two projections.

 Consider $\varphi \in L^1(\mu)$ and $\psi \in L^1(\nu)$, the Kantorovich dual problem is formulated as the following~\cite{Villani}:
\bq\label{eq:dual}
\sup_{\varphi, \psi} \left( \int_X \varphi\ d\mu + \int_Y \psi\ d\nu \right),
\eq
subject to $\varphi(x) + \psi(y) \leq c(x,y)$, for any $(x,y) \in X \times Y$.

The dual formulation is a linear optimization problem which is solvable by linear programming~\cite{cuturi2015smoothed,oberman2015efficient, schmitzer2016sparse}. Kantorovich obtained the 1975 Nobel prize in economics for his contributions to resource allocation problems where he interpreted the dual problem as an economic equilibrium. Recently Cuturi introduced the entropy regularized optimal transport problem which enforces the desirable properties for optimal transference plan and convexifies the problem. There have been extremely efficient computational algorithms~\cite{cuturi2013sinkhorn} which allow various applications in image processing, neuroscience, machine learning, etc~\cite{benamou2015iterative, cuturi2014fast, gramfort2015fast, solomon2015convolutional}.

In the 90's, Benamou and Brenier derived an equivalent dynamic formulation~\cite{BenBre} which has been one of the main tools for numerical computation. The Benamou-Brenier formula identifies the squared quadratic Wasserstein metric between $\mu$ and $\nu$ by 
\bq
W_2^2 (\mu, \nu) = \inf \int_0^1 \int |v(t,x)|^2 \rho(t,x) dx dt,
\eq
where the infimum is taken among all the solutions of the continuity equation:
\begin{eqnarray}
\frac{\partial \rho}{\partial t} + \nabla(v\rho) &=& 0, \\
\text{subject to\quad }\rho(0,x) = f,\  \rho(1,x) &=& g,\nonumber
\end{eqnarray}
In fact the infimum is taken among all Borel fields $v(t,x)$ that transports $\mu$ to $\nu$ continuously in time, satisfying the zero flux condition on the boundary. Many fast solvers based on this dynamic formulation has been proposed in literature~\cite{benamou2015augmented, li2017parallel, papadakis2014optimal}. They are used particularly in image registration, warping, texture mixing, etc. 

\subsubsection{The finite difference \MA solver}
As we have seen for the quadratic Wasserstein distance, the optimal map can be computed via the solution of a \MA partial differential equation~\cite{benamou2014numerical}. This approach has the advantage of drawing on the well-developed field of numerical partial differential equations (PDEs). We solve the \MA equation numerically for the viscosity solution using an almost-monotone finite difference method relying on the following reformulation of the \MA operator, which automatically enforces the convexity constraint~\cite{FroeseTransport}.
The scientific reason for using monotone type schemes follows from the following theorem by Barles and Souganidis~\cite{barles1991convergence}:
\begin{theorem}[Convergence of Approximation Schemes~\cite{barles1991convergence}]
Any consistent, stable, monotone approximation scheme  to the solution of fully nonlinear second-order elliptic or parabolic PDE converges uniformly on compact subsets to the unique viscosity solution of the limiting equation, provided this equation satisfies a  comparison principle.

\end{theorem}

The numerical scheme of~\cite{benamou2014numerical} uses the theory of~\cite{barles1991convergence} to construct a convergent discretization of the \MA equation~\eqref{eq:MAA} as stated in Theorem 7. A variational characterization of the determinant on the left hand side which also involves the negative part of the eigenvalues was proposed as the following equation:
\begin{multline}\label{eq:MA_convex}
{\det}(D^2u) = \\ \min\limits_{\{v_1,v_2\}\in V}\left\{\max\{u_{v_1,v_1},0\} \max\{u_{v_2,v_2},0\}+\min\{u_{v_1,v_1},0\} + \min\{u_{v_2,v_2},0\}\right\}
\end{multline}
where $V$ is the set of all orthonormal bases for $\R^2$.  

Equation~\eqref{eq:MA_convex} can be discretized by computing the minimum over finitely many directions $\{\nu_1,\nu_2\}$, which may require the use of a wide stencil.  
In the low-order version of the scheme, the minimum in~\eqref{eq:MA_convex} is approximated using only two possible values.  The first uses directions aligning with the grid axes.
\begin{multline}\label{MA1}
MA_1[u] = \max\left\{\Dt_{x_1x_1}u,\delta\right\}\max\left\{\Dt_{x_2x_2}u,\delta\right\} \\+ \min\left\{\Dt_{x_1x_1}u,\delta\right\} + \min\left\{\Dt_{x_2x_2}u,\delta\right\} - f / g\left(\Dt_{x_1}u, \Dt_{x_2}u\right) - u_0.
\end{multline}
Here $dx$ is the resolution of the grid, $\delta>K\Delta x/2$ is a small parameter that bounds second derivatives away from zero, $u_0$ is the solution value at a fixed point in the domain, and $K$ is the Lipschitz constant in the $y$-variable of $f(x)/g(y)$.

For the second value, we rotate the axes to align with the corner points in the stencil, which leads to
\begin{multline}\label{MA2}
MA_2[u] = \max\left\{\Dt_{vv}u,\delta\right\}\max\left\{\Dt_{\vp\vp}u,\delta\right\} + \min\left\{\Dt_{vv}u,\delta\right\} + \min\left\{\Dt_{\vp\vp}u,\delta\right\}\\ - f / g\left(\frac{1}{\sqrt{2}}(\Dt_{v}u+\Dt_{\vp}u), \frac{1}{\sqrt{2}}(\Dt_{v}u-\Dt_{\vp}u)\right) - u_0.
\end{multline}
Then the monotone approximation of the \MA equation is
\bq\label{eq:MA_compact} M_M[u] \equiv -\min\{MA_1[u],MA_2[u]\} = 0. \eq
We also define a second-order approximation, obtained from a standard centred difference discretisation,
\bq\label{eq:MA_nonmon} M_N[u] \equiv -\left((\Dt_{x_1x_1}u)(\Dt_{x_2x_2}u)-(\Dt_{x_1x_2}u^2)\right) + f/g\left(\Dt_{x_1}u,\Dt_{x_2}u\right) + u_0 = 0.\eq
These are combined into an almost-monotone approximation of the form
\bq\label{eq:MA_filtered} M_F[u] \equiv M_M[u] + \epsilon S\left(\frac{M_N[u]-M_M[u]}{\epsilon}\right) \eq
where $\epsilon$ is a small parameter and the filter $S$ is given by
\bq\label{eq:filter}
S(x) = \begin{cases}
x & \abs{x} \leq 1 \\
0 & \abs{x} \ge 2\\
-x+ 2  & 1\le x \le 2 \\
-x-2  & -2\le x\le -1.
\end{cases} 
\eq

The Neumann boundary condition is implemented using standard one-sided differences.  As described in~\cite{engquist2016optimal,FroeseTransport}, the (formal) Jacobian $\nabla M_F[u]$ of the scheme can be obtained exactly.  It is known to be sparse and diagonally dominant.


\begin{theorem}[Convergence to Viscosity Solution~{\cite[Theorem 4.4]{FroeseTransport}}]\label{thm:MA_convergence}
Let the \MA equation \eqref{eq:MAA} have a unique viscosity solution and let $g>0$ be Lipschitz continuous on $\R^d$. Then the solutions of the scheme \eqref{eq:MA_filtered} converge to the viscosity solution of \eqref{eq:MAA} with a formal discretization error of $\bO(Lh^2)$ where $L$ is the Lipschitz constant of $g$ and $h$ is the resolution of the grid.
\end{theorem} 

Once the discrete solution $u_h$ is computed, the squared Wasserstein metric is approximated via
\bq\label{eq:WassDiscrete}  
W_2^2(f,g) \approx \sum\limits_{j=1}^n (x_j-D_{x_j}u_h)^T\diag(f)(x_j-D_{x_j}u_h) dt,
\eq
where $n$ is the dimension of the data $f$ and $g$.
Then the gradient of the discrete squared Wasserstein metric can be expressed as
\bq 
\begin{split}
\frac{\partial W_2^2(f,g) }{\partial f} = &\sum\limits_{j=1}^n  \left[-2\nabla M_F^{-1}[u_f]^TD_{x_j}^T\diag(f) \right](x_j - D_{x_j}u_f)dt \\ & + \sum\limits_{j=1}^n  |x_j-D_{x_j}u_f|^2dt, 
\end{split}
\eq
This term is the discretized version of the Fr\'{e}chet derivative of the misfit function~\eqref{eqn:W22D} with respect to the synthetic data $f$, i.e., the adjoint source $\frac{\partial J}{\partial f}$ in the adjoint wave equation~\eqref{eq:FWI_adj}.




\section{Application of Optimal Transport to Seismic Inversion}
In this section, we first review the good properties of the $W_2$ norm for the application of full-waveform inversion. We will also explain some details of the implementations and show numerical results of using optimal transport based metrics as the misfit function in FWI.

\subsection{$W_2$ properties}

As we demonstrated in~\cite{engquist2016optimal}, the squared Wasserstein metric has several properties that make it attractive as a choice for misfit function.  One highly desirable feature is its convexity with respect to several parameterizations that occur naturally in seismic waveform inversion~\cite{yang2017application}. For example, variations in the  wave velocity lead to  simulated $f$ that are derived from shifts,
\bq\label{eq:shift}
f(x;s) = g(x+s\eta), \quad \eta \in \R^n,
\eq
 or dilations,
\bq\label{eq:dilation}
f(x;A) = g(Ax), \quad A^T = A, \, A > 0,
\eq
 applied to the observation $g$.
Variations in the strength of a reflecting surface or the focusing of seismic waves can also lead to local rescalings of the form
\bq\label{eq:rescaleLocal}
f(x;\beta) = \begin{cases}  \beta g(x), & x \in E\\ g(x), & x \in \R^n\backslash E.\end{cases}
\eq


\begin{theorem}[Convexity of squared Wasserstein metric~{\cite{engquist2016optimal}}]\label{thm:convexity}
The squared Wasserstein metric $W_2^2(f(m),g)$ is convex with respect to the model parameters $m$ corresponding to a shift~$s$ in~\eqref{eq:shift}, the eigenvalues of a dilation matrix~$A$ in~\eqref{eq:dilation}, or the local rescaling parameter~$\beta$ in~\eqref{eq:rescaleLocal}.
\end{theorem} 

Another important property of optimal transport is the insensitivity to noise. All seismic data contains either natural or experimental equipment noise. For example, the ocean waves lead to extremely low-frequency data in the marine acquisition. Wind and cable motions also generate random noise.

\begin{theorem}[Insensitivity to noise~{\cite{engquist2016optimal}}]\label{thm:noise}
Let $f_{ns}$ be $f$ with a piecewise constant additive noise of mean zero uniform distribution.
The squared Wasserstein metric $W_2^2(f,f_{ns})$ is of $\bO(\frac{1}{N})$ where $N$ is the number of pieces of the additive noise in $f_{ns}$.
\end{theorem} 

The $L^2$ norm is known to be sensitive to noise since the misfit between clean and noisy data is calculated as the sum of squared noise amplitude at each sampling point. 


\subsection{Data normalization}
In optimal transport theory, there are two main requirements for signals $f$ and $g$: positivity and mass balance.  Since these are not expected for seismic signals, some data pre-processing is needed before we can implement Wasserstein-based FWI.
In~\cite{EFWass,engquist2016optimal}, the signals were separated into positive and negative parts $f^+ = \max\{f,0\}$, $f^- = \max\{-f,0\}$ and scaled by the total mass $\langle f \rangle = \int_X f(x)\,dx$.  Inversion was accomplished using the modified misfit function
\bq W_2^2\left(\frac{f^+}{\langle f^+ \rangle},\frac{g^+}{\langle g^+ \rangle} \right) + W_2^2\left(\frac{f^-}{\langle f^- \rangle}, \frac{g^-}{\langle g^- \rangle}\right). \eq

While this approach preserves the desirable theoretical properties of convexity to shifts and noise insensitivity, it is not easy to combine with the adjoint-state method and more realistic examples. We require the scaling function to be differentiable so that it is easy to apply the chain rule when calculating the Fr\'{e}chet derivative for FWI backpropagation and also better suited for the \MA and the wave equation solvers. 

There are other ways to rescale the datasets so that they become positive. For example, we can square the data as $\tilde{f} = f^2$ or extract the envelope of the data. These methods preserve the convexity concerning simple shifts, but we have lost the uniqueness: $f^2 = g^2$ does not imply $f=g$. As a result, more local minima are present since the fact that the misfit $J(f^2, g^2)$ is decreasing does not necessarily indicate that $f$ is approaching $g$, not to mention the non-unique issue of the inverse problem itself.

\begin{figure}
\centering
  \subfloat[Linear]{\includegraphics[width=0.3\textwidth]{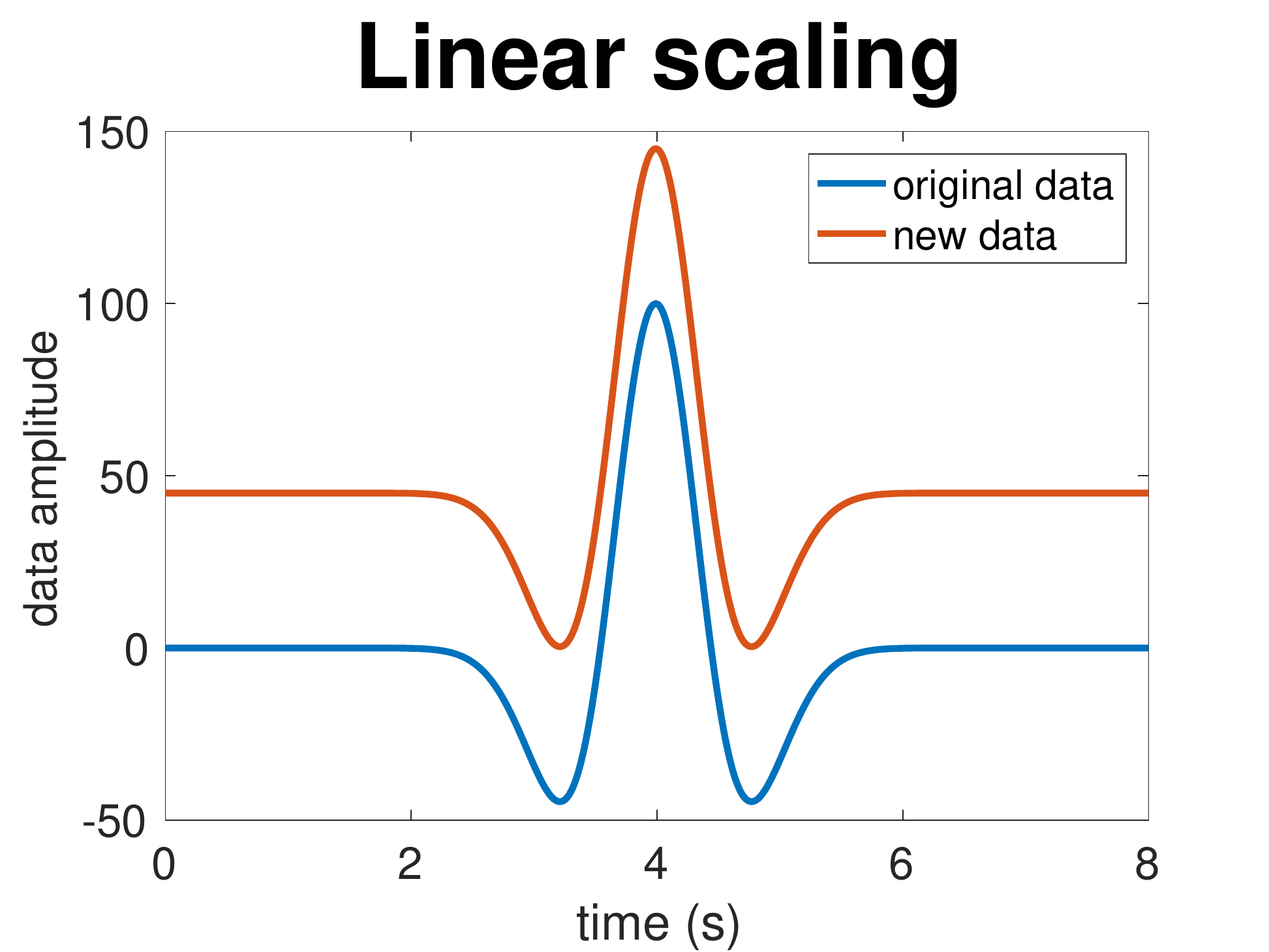}\label{fig:norm_LT}}
  \subfloat[Exponential]{\includegraphics[width=0.3\textwidth]{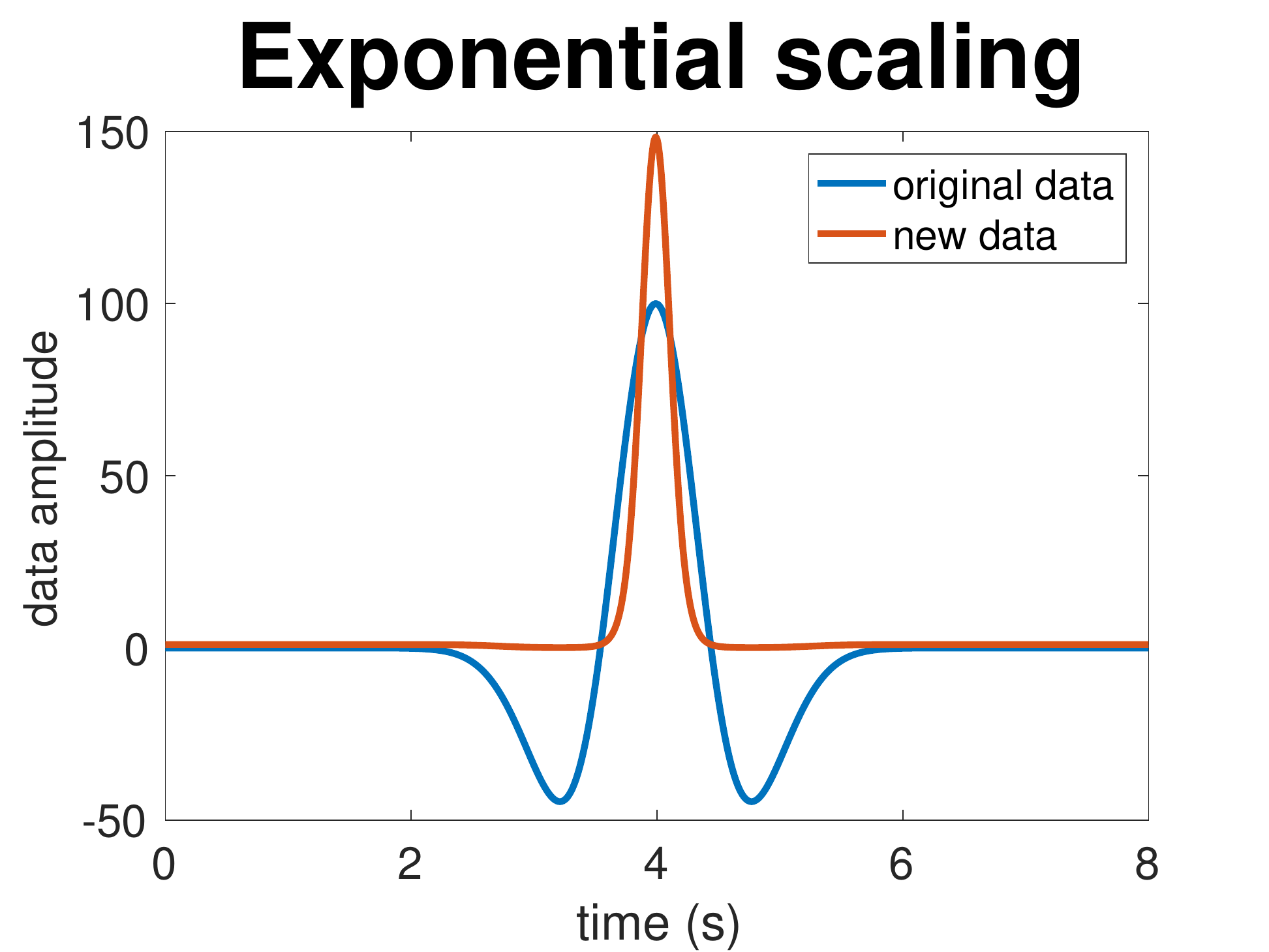}\label{fig:norm_exp}}
  \subfloat[Sign-sensitive]{\includegraphics[width=0.3\textwidth]{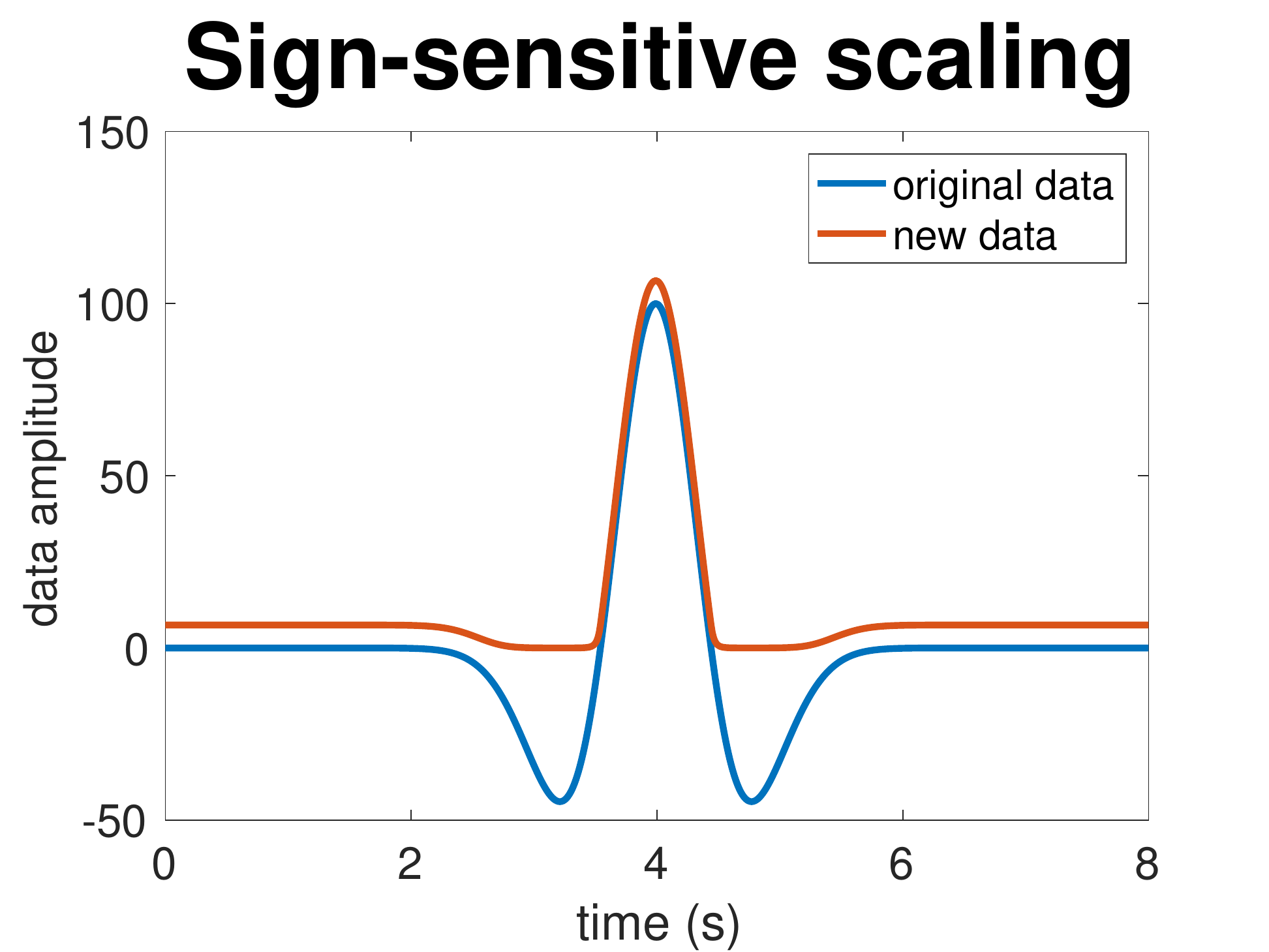}\label{fig:norm_mix}}
  \caption{(a)~The linear, (b)~the exponential and (c) the sign-sensitive scaling of a Ricker wavelet (Blue)}
  \label{fig:data_norm}
\end{figure}

Typically, we first scale the data $f$ to be positive as $\tilde{f}$ and then normalize to ensure mass balance as $\tilde{f} / <\tilde{f}>$. We now introduce three normalization methods that are robust in realistic large-scale inversions: the linear scaling~\cite{yang2017application} (Figure ~\ref{fig:norm_LT})
\bq\label{eq:norm_LT}
\tilde{f} = f  + c_1,\quad c_1 \geq \max\{-f,-g\},
\eq
the exponential scaling~\cite{qiu2017full} (Figure ~\ref{fig:norm_exp})
\bq\label{eq:norm_exp}
\tilde{f} =\exp(c_2 f),\quad c_2>0,
\eq
and the sign-sensitive scaling (Figure ~\ref{fig:norm_mix})
\bq\label{eq:norm_mix}
    \tilde{f}=
             \begin{cases}
                  f + \frac{1}{c_3},\quad & f\geq 0  \\
                \frac{1}{c_3} \exp(c_3f),\quad & f < 0 \\
            \end{cases},  \quad c_3 >0.
\eq

If $c_2$ in ~\eqref{eq:norm_exp} and $c_3$ in ~\eqref{eq:norm_mix} are large enough, these two scaling methods keep the convexity of $W_2$ norm regarding simple shifts as shown in Figure~\ref{fig:2_ricker_W2}.  From Taylor expansion, we can see that the scalings are very close to the linear scaling when $c_2$ is small. One has to be careful with the exponential scaling~\eqref{eq:norm_exp} since it can easily become extremely large, but the sign-sensitive scaling~\eqref{eq:norm_mix} will not.

\subsection{FWI with Kantorovich-Rubinstein norm}
When the cost function $c(x,y)$ is the $L^1$ norm $|x-y|$, i.e. $p=1$ in ~\eqref{eq:static} with $f\geq 0$, $g\geq 0$, and $\int f = \int g$, the corresponding alternative $W_1$ distance has the following equivalent dual formulation:
\bq \label{eq:W1}
W_1(f,g)  
=\max_{\varphi \in \text{Lip}_1} \int_X \varphi(x)(f(x) - g(x))dx, 
\eq
where $\text{Lip}_1$ is the space of Lipschitz continuous functions with Lipschitz constant 1. However, seismic data $f$ and $g$ are oscillatory containing both positive and negative parts. If $\int f \neq \int g$, the value of ~\eqref{eq:W1} is always $+\infty$. Recently, \cite{W1_2D, W1_3D} introduced the following Kantorovich-Rubinstein (KR) norm in FWI which is a relaxation of the original $W_1$ distance by constraining the dual space:
\begin{equation}  \label{eq: KR}
\text{KR}(f,g) = \max_{\varphi \in \text{BLip}_1} \int_X \varphi(x)(f(x) - g(x))dx
\end{equation}
Here $\text{BLip}_1$ is the space of bounded Lipschitz continuous functions with Lipschitz constant 1.
One advantage of using KR norm in FWI is that there is no need to normalize the data to be positive and mass balanced. However, KR norm has no direct connection with optimal transport once we no longer require $f$ and $g$ to be probability measures~\cite{vershik2013long}. When $f$ and $g$ are far apart which is very common when the initial velocity is rough, the maximum in ~\eqref{eq: KR} is achieved by ``moving'' $f^+$ to $f^-$ and $g^+$ to  $g^-$. The notion of transport is void in this case and convexity is lost.

\subsection{Numerical results of global $W_2$}
In the next two subsections, we provide numerical results for two approaches to using $W_2$ with linear normalization~\eqref{eq:norm_LT}: trace-by-trace comparison and using the entire 2D datasets as objects.  Here a trace is the time history measured at one receiver while the entire dataset consists of the time history of all the receivers. These are compared with results produced by using the standard least-squares norm $L^2$ to measure the misfit. More examples can be found in~\cite{yang2017application}.

\begin{figure}
\centering
  \subfloat[]{\includegraphics[width=0.45\textwidth]{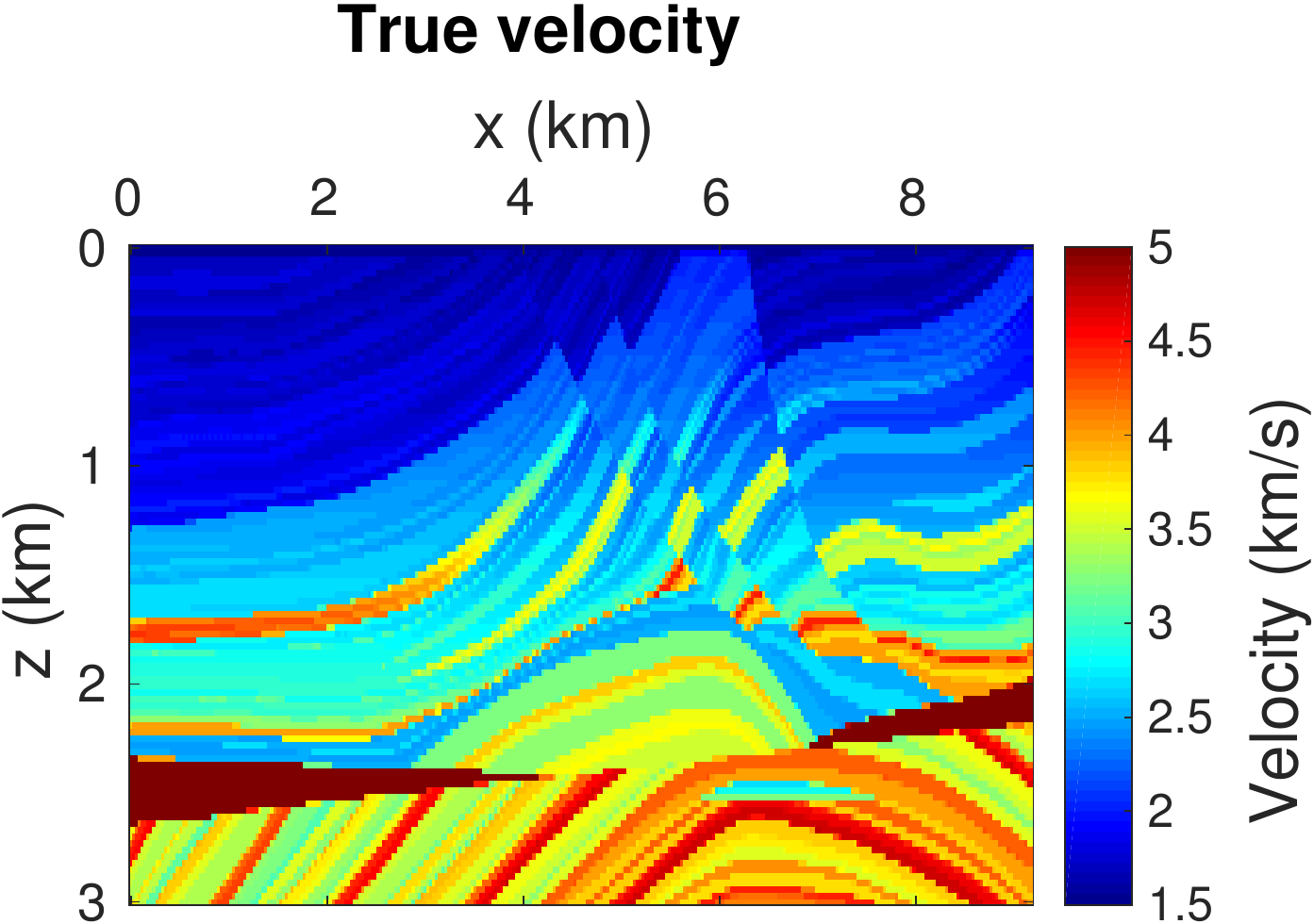}\label{fig:marm2_true}}
  \subfloat[]{\includegraphics[width=0.45\textwidth]{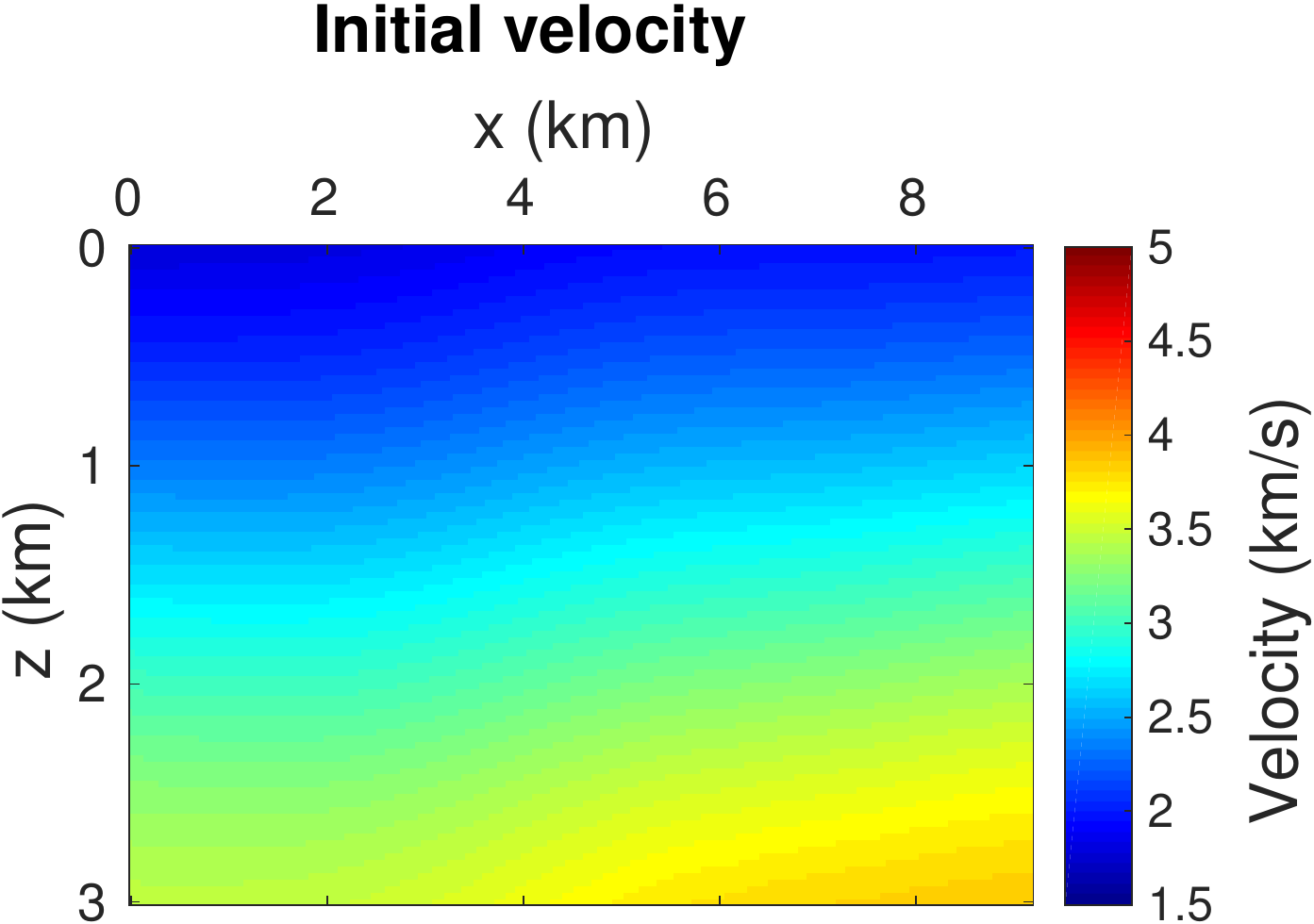}\label{fig:marm2_v0}}
  \caption{(a)~True velocity and (b)~inital velocity for full Marmousi model}
  \label{fig:marm2_true,marm2_v0}
\end{figure}


\begin{figure}
\centering
  \subfloat[]{\includegraphics[width=0.45\textwidth]{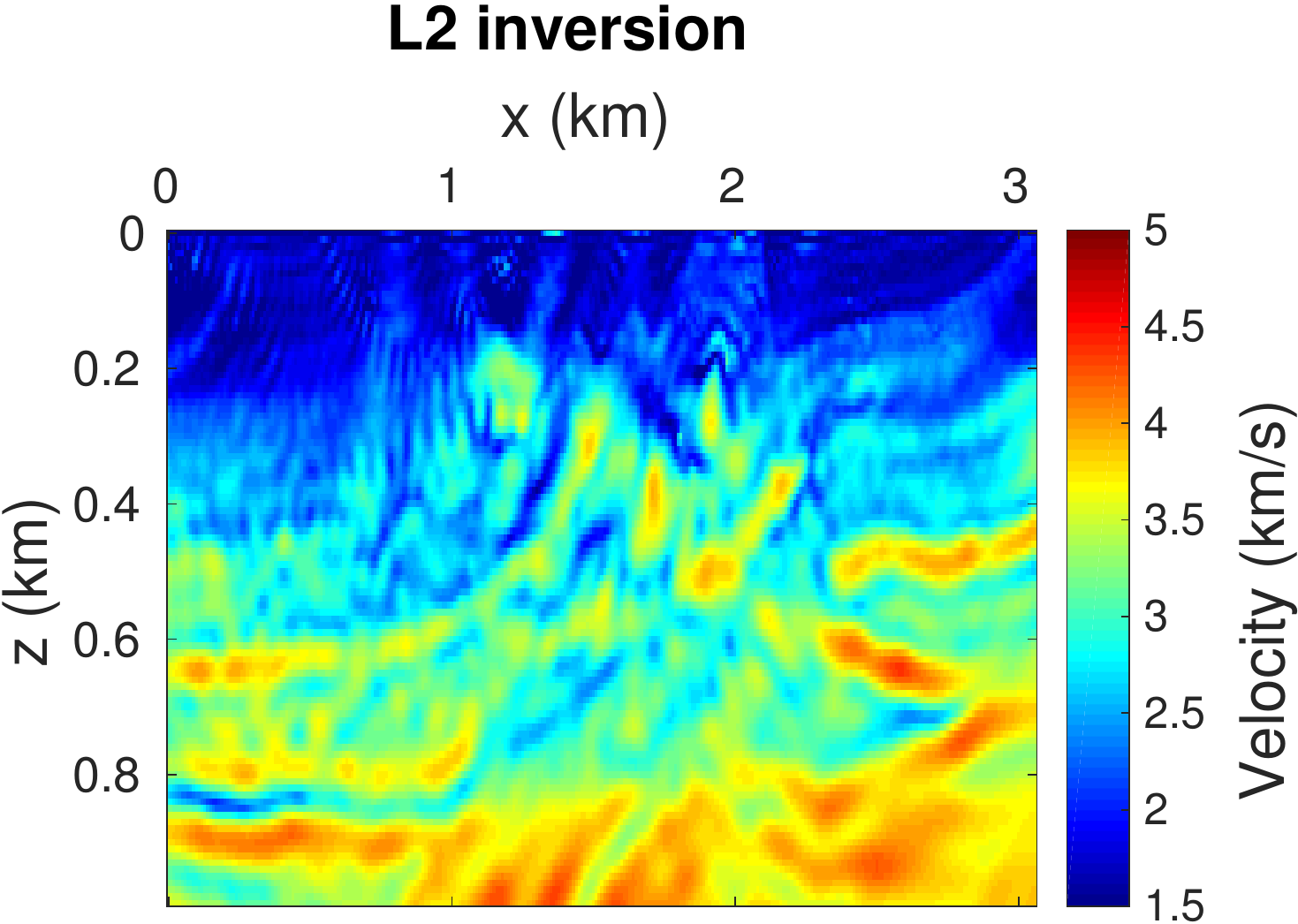}\label{fig:marm_L2}}
  \subfloat[]{\includegraphics[width=0.45\textwidth]{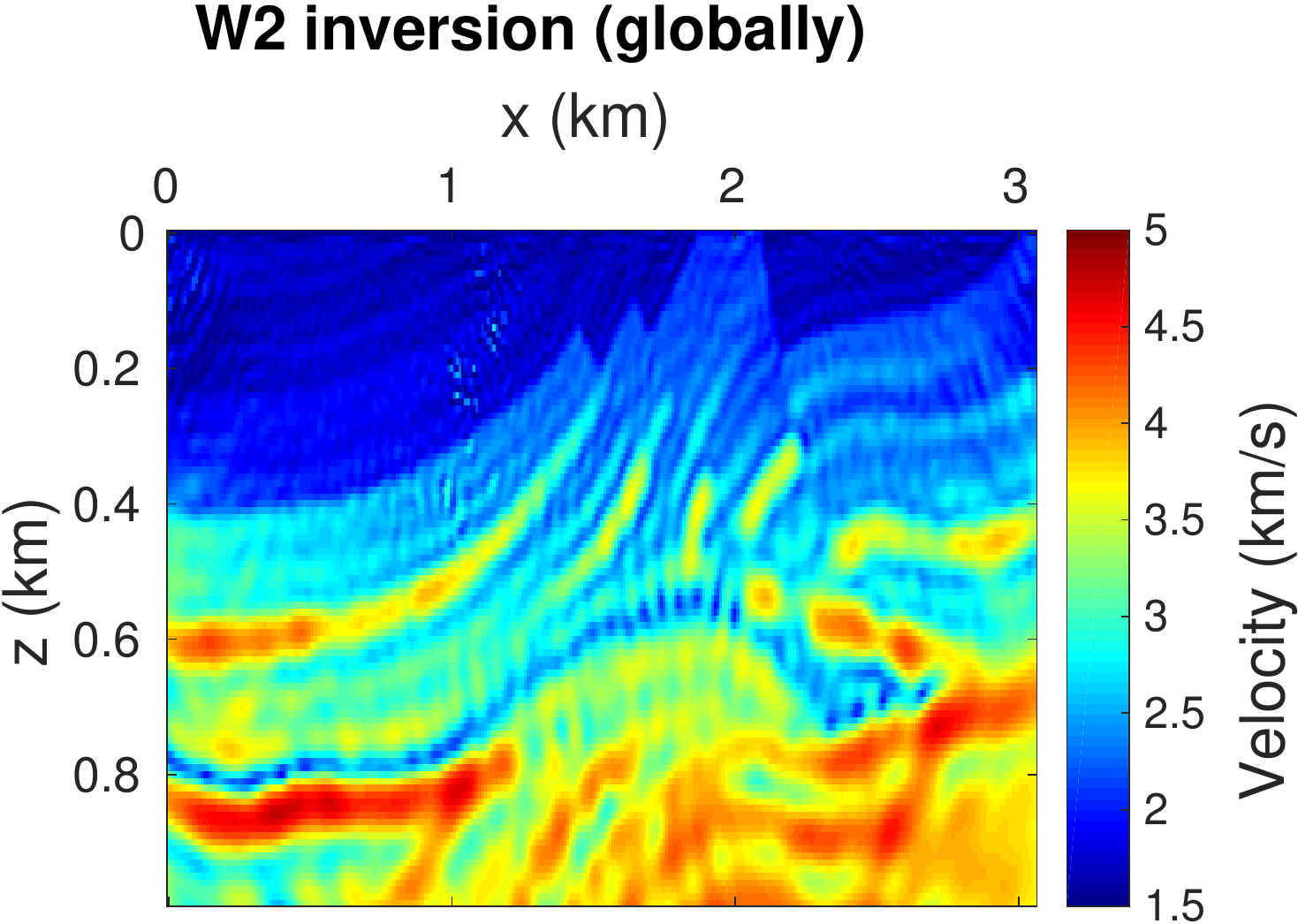}\label{fig:marm_w2_2D}}
  \caption{Inversion results of (a)~$L^2$ and (b)~global $W_2$  for the scaled Marmousi model}
  \label{fig:marm_inv_scaled}
\end{figure}

First, we use a scaled Marmousi model to compare the inversion between global $W_2$ and the conventional $L^2$ misfit function. Figure~\ref{fig:marm2_true} is the P-wave velocity of the true Marmousi model, but in this experiment, we use a scaled model which is 1 km in depth and 3 km in width. The inversion starts from an initial model that is the true velocity smoothed by a Gaussian filter with a deviation of 40, which is highly smoothed and far from the true model (a scaled version of Figure~\ref{fig:marm2_v0}). We place 11 evenly spaced sources on top at 50 m depth and 307 receivers on top at the same depth with a 10 m fixed acquisition. The discretization of the forward wave equation is 10 m in the $x$ and $z$ directions and 10 ms in time. The source is a Ricker wavelet which is the second derivative of the Gaussian function with a peak frequency of 15 Hz, and a bandpass filter is applied to remove the frequency components from 0 to 2 Hz. 


We use L-BFGS, a quasi-Newton method as the optimization algorithm ~\cite{liu1989limited}.
Inversions are terminated after 200 iterations. Figure~\ref{fig:marm_L2} shows the inversion result using the traditional $L^2$ least-squares method after 200 L-BFGS iterations. The inversion result of global $W_2$ (Figure~\ref{fig:marm_w2_2D}) avoids the problem of local minima suffered by the conventional $L^2$ metric, whose result demonstrates spurious high-frequency artifacts due to a point-by-point comparison of amplitude. 

We solve the \MA equation numerically in each iteration of the inversion. The drawback to the PDE approach is that data must be sufficiently regular for solutions to be well-defined and for the numerical approximation to be accurate.  To remain robust on realistic examples, we use filters that effectively smooth the seismic data, which can lead to a loss of high-frequency information.  For illustration in this paper, we perform computations using a \MA solver for synthetic examples.  Even in 2D, some limitations are apparent. This is expected to become even more of a problem in higher-dimensions and motivates our introduction of a trace-by-trace technique that relies on the exact 1D solution. The trace-by-trace technique is currently more promising for practical applications, as is evidenced in our computational examples in the next section.

\subsection{Numerical results of trace-by-trace $W_2$}
Recall that for the 1D trace-by-trace approach, the misfit function in FWI is
\bq 
J_1(m) = \sum\limits_{r=1}^R W_2^2(f(\mathbf{x_r},t;m),g(\mathbf{x_r},t)), 
\eq
where $R$ is the total number of traces, $g$ is observed data, $f$ is simulated data, $\mathbf{x_r}$ are receiver locations, and $m$ is the model parameter. The adjoint source term for each single trace is 
\bq~
\frac{\partial W^2_2(f,g)}{\partial f} = \left( \int_t^{T_0}\frac{-2(s-G^{-1} \circ F(s))}{g(G^{-1}\circ F(s))} f(s) ds+ |t-G^{-1}(F(t))|^2 \right)dt.
\eq

\begin{figure}
\centering
  \subfloat[]{\includegraphics[width=0.45\textwidth]{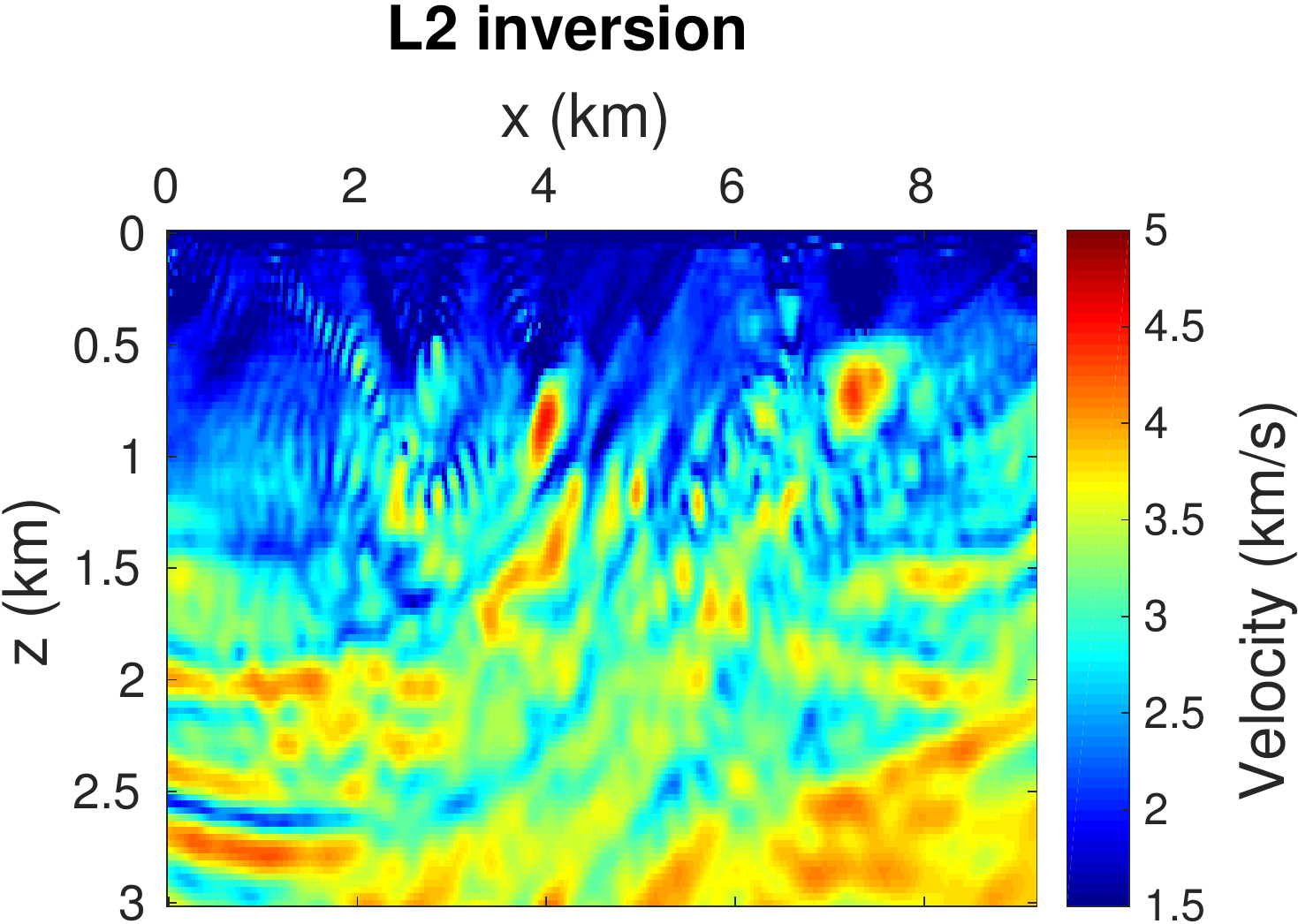}\label{fig:marm2_L2}}
  \subfloat[]{\includegraphics[width=0.45\textwidth]{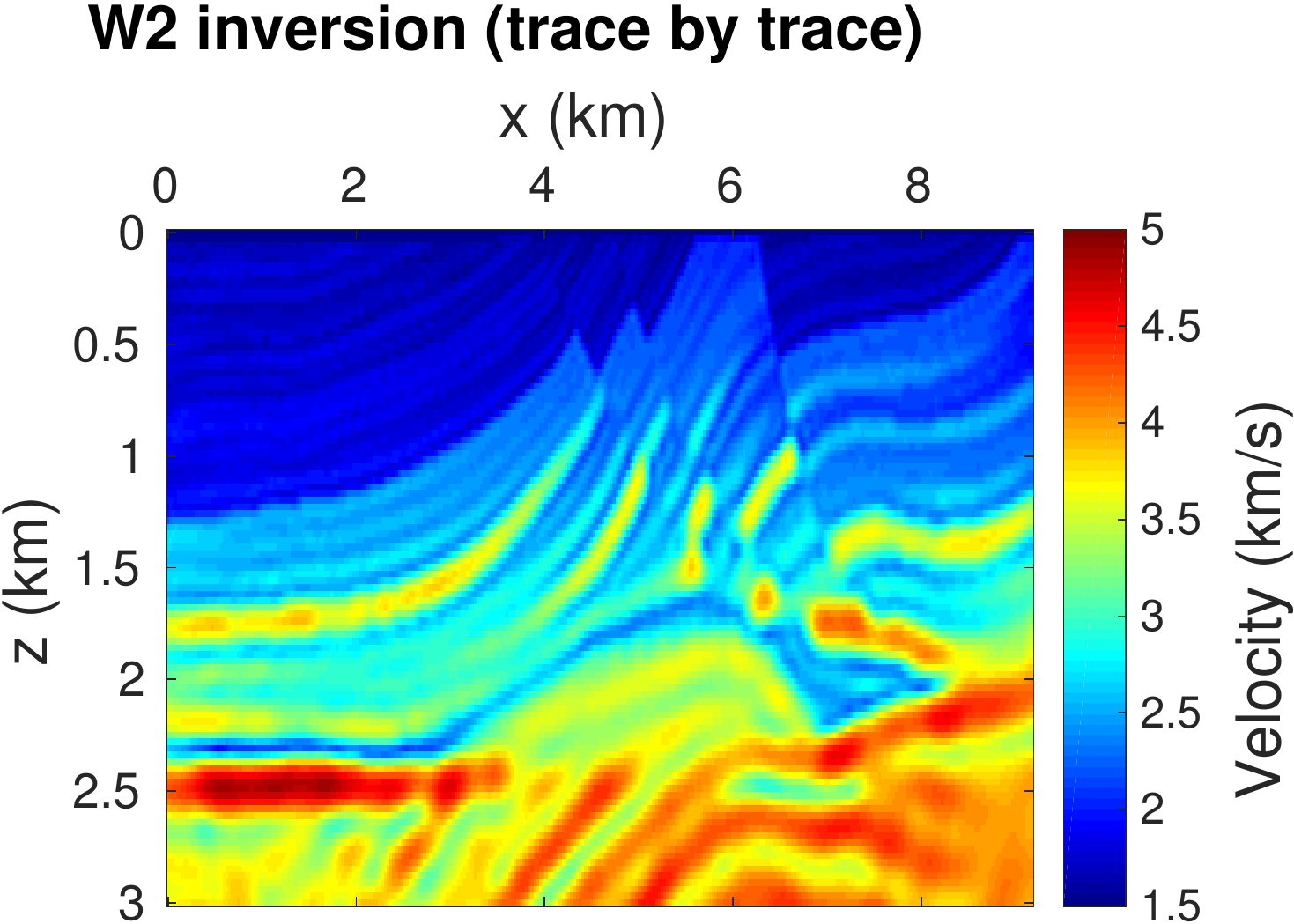}\label{fig:marm2_w2_1D}}
  \caption{Inversion results of (a)~$L^2$ and (b)~trace-by-trace $W_2$ for the true Marmousi model}
  \label{fig:marm2_inv}
\end{figure}

\begin{figure}
\centering
 {\includegraphics[width=1.0\textwidth]{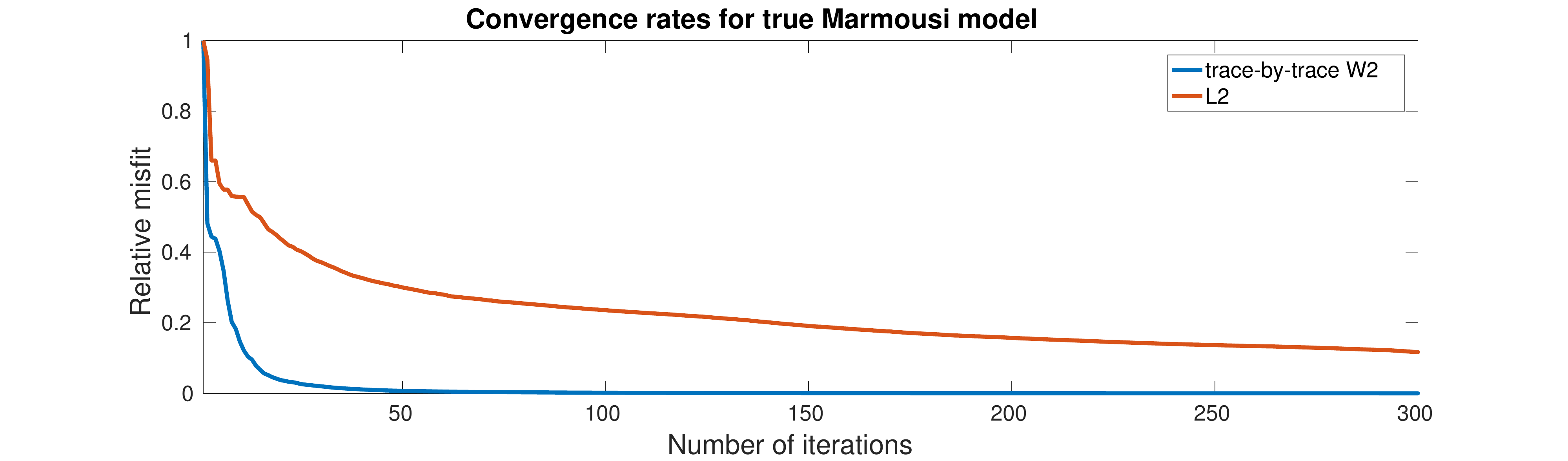}}
   \caption{The convergence curves for trace-by-trace $W_2$ and $L^2$ based inversion of the full Marmousi model}
    \label{fig:marm2_conv}
\end{figure}

The next experiment is to invert the full Marmousi model by conventional $L^2$ and trace-by-trace $W_2$ misfit. Figure~\ref{fig:marm2_true} is the P-wave velocity of the true Marmousi model, which is 3 km in depth and 9 km in width. The inversion starts from an initial model that is the true velocity smoothed by a Gaussian filter with a deviation of 40 (Figure~\ref{fig:marm2_v0}). The rest of the settings are the same as the previous section. Inversions are terminated after 300 L-BFGS iterations. Figure~\ref{fig:marm2_L2} shows the inversion result using the traditional $L^2$ least-squares method and figure~\ref{fig:marm2_w2_1D} shows the final result using trace-by-trace $W_2$ misfit function. Again, the result of $L^2$ metric has spurious high-frequency artifacts while $W_2$ correctly inverts most details in the true model. The convergence curves in Figure~\ref{fig:marm2_conv} show that $W_2$ reduces the relative misfit to 0.1 in 20 iterations while $L^2$ converges slowly to a local minimum.

\subsection{Insensitivity to noise}
\begin{figure}
\centering
  \subfloat[]{\includegraphics[height=4cm]{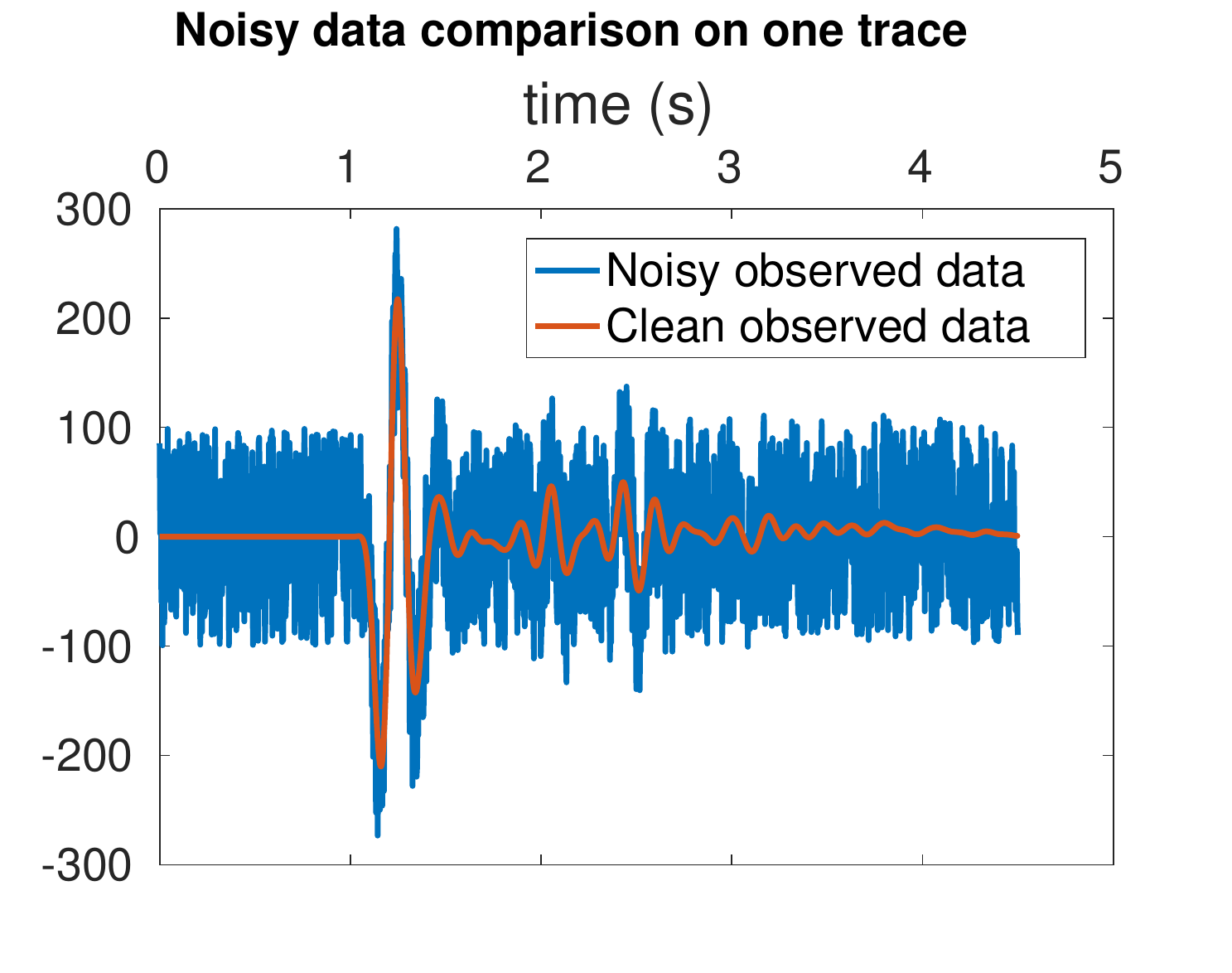}\label{fig:MARM2_noisy_trace}}
  \subfloat[]{\includegraphics[height=5cm]{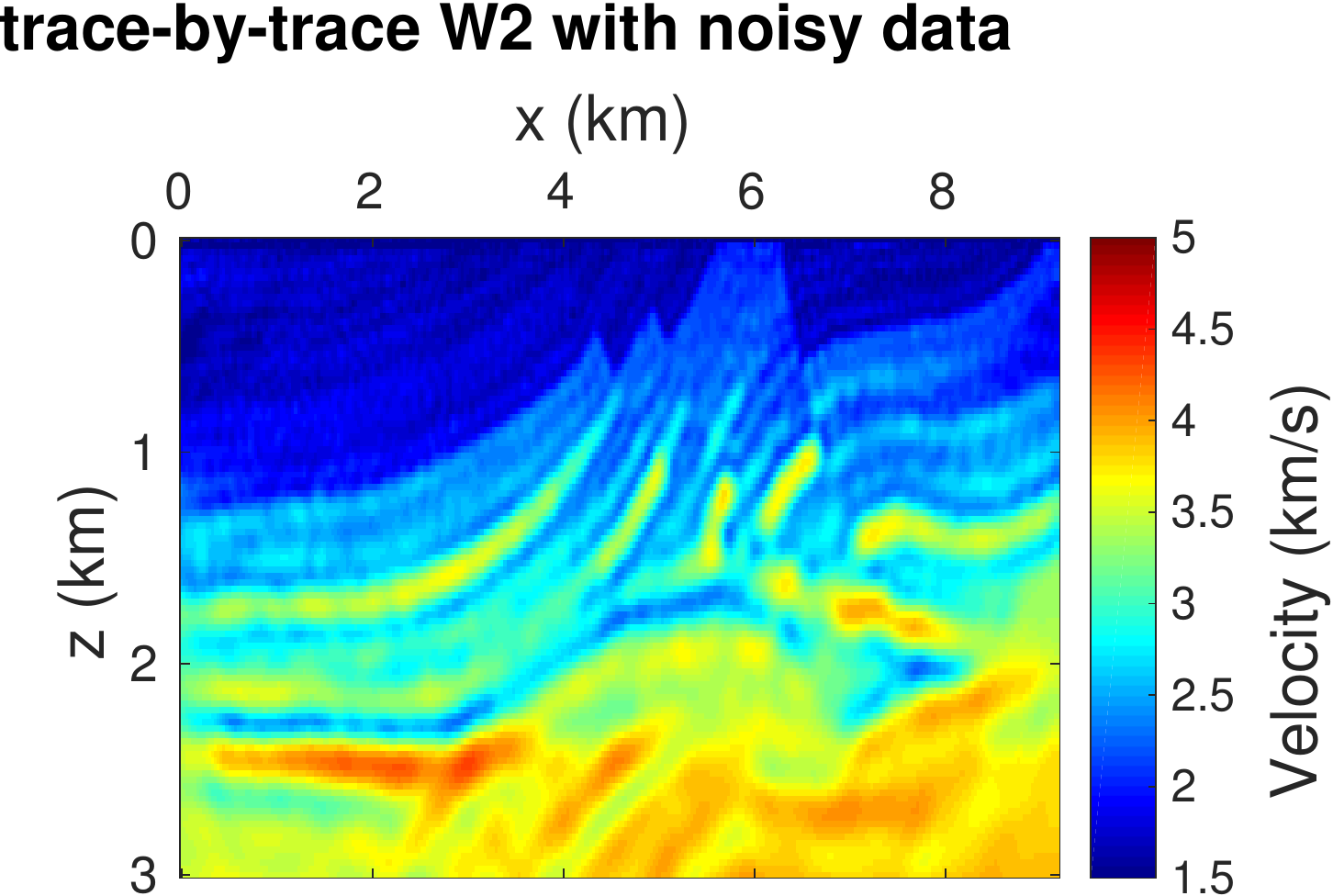}\label{fig:MARM2_noisy_vel}}
  \caption{(a)~Noisy and clean data and (b)~inversion result with the noisy data}
  \label{fig:marm2_noisy}
\end{figure}

One of the good properties of the \QW is the insensitivity to noise~\cite{engquist2016optimal}. We repeat the previous experiment with a noisy reference by adding a uniform random iid noise to the data from the true velocity (Figure~\ref{fig:MARM2_noisy_trace}). The signal-to-noise ratio (SNR) is $-3.47$ dB. In optimal transport, the effect of noise is essentially negligible due to the strong cancellation between the nearby positive and negative noise.

All the settings remain the same as in the previous numerical experiment except the observed data. After 96 iterations, the optimization converges to a velocity presented in Figure~\ref{fig:MARM2_noisy_vel}. Although the result has lower resolution than Figure~\ref{fig:marm2_w2_1D}, it still recovers most features of Marmousi model correctly. Even when the noise is much larger than the signal, the \QW still converges reasonably well.

\bibliographystyle{plain}

\end{document}